\documentclass[a4paper, 10pt]{article}
\usepackage{cancel}
\usepackage[normalem]{ulem}

    \makeatletter
      
\usepackage{fancyhdr}
\usepackage{bm}
\usepackage{cancel}
\usepackage{latexsym}
    
\usepackage{amssymb,amsmath,amsfonts,bbm,pifont,upgreek,bbold,accents}  
\usepackage[colorlinks=true]{hyperref}
\hypersetup{urlcolor=blue, citecolor=red}

\usepackage[dvips]{graphicx}
\font\teneufm=eufm10
\font\seveneufm=eufm7
\font\fiveeufm=eufm5
\newfam\eufmfam
\textfont\eufmfam=\teneufm
\scriptfont\eufmfam=\seveneufm
\scriptscriptfont\eufmfam=\fiveeufm

\newcommand\beq[1]{\begin{equation}\label{#1} }
\newcommand{\eeq}{\end{equation} }

\newcommand\beqa[1]{\begin{eqnarray} \label{#1}}
\newcommand{\eeqa}{\end{eqnarray} }
\newcommand{\beqano}{\begin{eqnarray*} }
\newcommand{\eeqano}{\end{eqnarray*} }
\newcommand\arr[1]{\left\{\begin{array}{l}#1\end{array}\right.}
\renewcommand{\theequation}{\arabic{section}.\arabic{equation}}

\newtheorem{theorem}{Theorem}[section]
\newtheorem{definition}{Definition}[section]
\newtheorem{proposition}{Proposition}[section]
\newtheorem{lemma}{Lemma}[section]

\newtheorem{sublemma}{Sublemma}[section]
\newtheorem{remark}{Remark}[section]
\newtheorem{notationalremark}{Notational Remark}[section]
\newtheorem{corollary}{Corollary}[section]
\newtheorem{assumption}{Assumption}[section]
\newtheorem{claim}{Claim}[section]

\newtheorem{tools}{$\negsp\negsp$}[subsection]

\newcommand\thm[1]{\begin{theorem}\label{#1}}
\newcommand\thmtwo[2]{\begin{theorem}[#1]\label{#2}}
\newcommand\ethm{\end{theorem} }
\newcommand\dfn[1]{\begin{definition}\label{#1} \rm}
\newcommand\dfntwo[2]{\begin{definition}[#1]\label{#2} \rm}
\newcommand\edfn{\end{definition} }
\newcommand\pro[1]{\begin{proposition}\label{#1}}
\newcommand\protwo[2]{\begin{proposition}[#1]\label{#2}}
\newcommand\epro{\end{proposition} }
\newcommand\lem[1]{\begin{lemma}\label{#1}}
\newcommand\lemtwo[2]{\begin{lemma}[#1]\label{#2}}
\newcommand\elem{\end{lemma} }
\newcommand\sublem[1]{\begin{sublemma}\label{#1}}
\newcommand\sublemtwo[2]{\begin{sublemma}[#1]\label{#2}}
\newcommand\esublem{\end{sublemma} }
\newcommand\rem[1]{\begin{remark}\label{#1} \rm}
\newcommand\erem{\end{remark} }
\newcommand\notrem[1]{\begin{notationalremark}\label{#1} \rm}
\newcommand\enotrem{\end{notationalremark} }
\newcommand\cor[1]{\begin{corollary}\label{#1}}
\newcommand\cortwo[2]{\begin{corollary}[#1]\label{#2}}
\newcommand\ecor{\end{corollary} }
\newcommand\asmp[1]{\begin{assumption}\label{#1}}
\newcommand\asmptwo[2]{\begin{assumption}[#1]\label{#2}}
\newcommand\easmp{\end{assumption} }
\newcommand\clm[1]{\begin{claim}\label{#1}}
\newcommand\eclm{\end{claim} }

\newcommand\equ[1]{{\rm (\ref{#1})}}

\expandafter\chardef\csname pre amssym.def
at\endcsname=\the\catcode`\@
\catcode`\@=11
\def\undefine#1{\let#1\undefined}
\def\newsymbol#1#2#3#4#5{\let\next@\relax
 \ifnum#2=\@ne\let\next@\msafam@\else
 \ifnum#2=\tw@\let\next@\msbfam@\fi\fi
 \mathchardef#1="#3\next@#4#5}
\def\mathhexbox@#1#2#3{\relax
 \ifmmode\mathpalette{}{\m@th\mathchar"#1#2#3}%
 \else\leavevmode\hbox{$\m@th\mathchar"#1#2#3$}\fi}
\def\hexnumber@#1{\ifcase#1 0\or 1\or 2\or 3\or 4\or 5\or 6\or 7\or
8\or
 9\or A\or B\or C\or D\or E\or F\fi}
\ifcase\@ptsize
 \font\tenmsb=msbm10
 \font\sevenmsb=msbm7
 \font\fivemsb=msbm5
\or
 \font\tenmsb=msbm10 scaled \magstephalf
 \font\sevenmsb=msbm7 scaled \magstephalf
 \font\fivemsb=msbm5  scaled \magstephalf
\or
 \font\tenmsb=msbm10 scaled \magstep1
 \font\sevenmsb=msbm7 scaled \magstep1
 \font\fivemsb=msbm5 scaled \magstep1
\fi
\newfam\msbfam
\textfont\msbfam=\tenmsb
\scriptfont\msbfam=\sevenmsb
\scriptscriptfont\msbfam=\fivemsb
\edef\msbfam@{\hexnumber@\msbfam}
\def\Bbb#1{\fam\msbfam\relax#1}
\def\widehat#1{\setboxz@h{$\m@th#1$}%
 \ifdim\wdz@>\tw@ em\mathaccent"0\msbfam@5B{#1}%
 \else\mathaccent"0362{#1}\fi}
\def\widetilde#1{\setboxz@h{$\m@th#1$}%
 \ifdim\wdz@>\tw@ em\mathaccent"0\msbfam@5D{#1}%
 \else\mathaccent"0365{#1}\fi}

\def\RIfM@{\relax\ifmmode}
\def\nonmatherr@#1{\errmessage{\string#1\space allowed only in math mode}}
\def\Bbb{\RIfM@\expandafter\Bbb@\else
 \expandafter\nonmatherr@\expandafter\Bbb\fi}
\def\Bbb@#1{{\Bbb@@{#1}}}
\def\Bbb@@#1{\fam\msbfam\relax#1}
\def\setboxz@h{\setbox\z@\hbox}
\def\wdz@{\wd\z@}
\catcode`\@=\csname pre amssym.def at\endcsname
\newcommand{\nl}{{\smallskip\noindent}}

\newcommand{\negsp}{\hspace{-.09truecm}}  

\newcommand{\dst}{\displaystyle}
\newcommand\ovl[1]{\overline {#1} }

\newcommand\su[1]{\frac{1}{{#1}} }

\newcommand\sign{{\, \rm sign\, }}

\newcommand{\tg}{{\rm \, tg \, }}
\newcommand{\torus}{{\Bbb T}   }

\newcommand{\real}{{\Bbb R}   }
\newcommand{\integer}{{\Bbb Z}   }
\newcommand{\complex}{{\Bbb C}   }

\renewcommand{\a }{{\alpha}   }
\renewcommand{\b}{{\beta}   }
\newcommand{\g}{{\gamma}   }
\newcommand{\G}{{\Gamma}   }
\renewcommand{\d}{{\delta}   }

\renewcommand{\l}{{\lambda}   }
\renewcommand{\L}{{\Lambda}   }

\newcommand{\n}{{\nu}   }

\newcommand{\p}{{\pi}   }

\newcommand{\s}{{\sigma}   }

\renewcommand{\t}{{\tau}   }
\newcommand{\f}{{\varphi}   }

\renewcommand{\o}{{\omega}   }

\newcommand{\const}{{\, \rm const}}

\newcommand{\cA}{{\cal A} }
\newcommand{\cB}{{\cal B} }
\newcommand{{\cE}}{{\cal  E} }
\newcommand{\cT}{{\cal T} }
\newcommand{\cR}{{\cal R} }
\newcommand{{\cH}}{{\cal H} }
\newcommand{{\cK}}{{\cal K} }
\newcommand{\cC}{{\cal C} }
\newcommand{\cD}{{\cal D} }
\newcommand{\cF }{{\cal F} }
\newcommand{\cG}{{\cal G} }
\newcommand{{\cJ}}{{\cal J}}
\newcommand{\cL}{{\cal L} }
\newcommand{\cM}{{\cal M} }
\newcommand{\cP}{{\cal P} }
\newcommand{\cI}{{\cal I} }

\newcommand{\cS}{{\cal S} }

\newcommand\bx{{\mathbf x}}

\newcommand\by{{\mathbf y}}

\newcommand\ppu{{(1) }}
\newcommand\ppd{{(2) }}

\newcommand{\CC}{{\rm C}}

\newcommand{\EE}{{\rm E}}
\newcommand\FF{{\rm F}}
\newcommand\GG{{\rm G}}
\newcommand\HH{{\rm H}}
\newcommand{\II}{{\rm I}}
\newcommand\JJ{{\rm J}}
\newcommand\KK{{\rm K}}

\newcommand\MM{{\rm M}}

\newcommand\OO{{\rm O}}

\newcommand\RR{{\rm R}}

\newcommand\UU{{\rm U}}

\newcommand\XX{{\rm X}}
\newcommand\YY{{\rm Y}}
\newcommand\ZZ{{\rm Z}}

\newcommand\ee{{\rm e}}
\newcommand\hh{{\rm h}}

\newcommand\rr{{\rm r}}

\newcommand\zz{{\rm z}}
\usepackage{color}
\definecolor{yellow}{rgb}{0.99, 0.93, 0.0}

\newcommand\mm{{\rm m}}

\begin{document}

\title{\bf 
Euler integral and perihelion librations\thanks{The  author is supported by  the European Research Council. Grant 677793 Stable and Chaotic Motions in the Planetary Problem
{\bf MSC2000 numbers:}
primary:
34C20, 70F07,  37J10, 37J15, 37J35;
secondary: 
34D10,  70F10, 70F15, 37J25, 37J40. {\bf Keywords:} Two--centre problem; Euler Integral; Canonical coordinates; three--body problem.}}

\author{ 
Gabriella Pinzari
\\
\vspace{-.2truecm}
\\{\footnotesize Dipartimento di Matematica T. Levi--Civita}
\vspace{-.2truecm}
\\{\footnotesize via Trieste, 63, 35131, Padova (Italy)}
\vspace{-.2truecm}
\\{\scriptsize gabriella.pinzari@math.unipd.it}
}\date{February,  2020}
\maketitle

\begin{abstract}\footnotesize{We discuss  dynamical aspects  of an analysis of the two--centre problem  started in~\cite{pinzari19}. The perturbative nature of our approach allows us to foresee  applications to the three--body problem.}

\end{abstract}

\maketitle

\tableofcontents


\renewcommand{\theequation}{\arabic{equation}}
\setcounter{equation}{0}

\newpage

\section{Introduction}\label{2centres}

The two--centres (or {\it Euler--}) problem is the $3$--degrees of freedom  ($2$ in the plane) system of one particle interacting with two fixed masses via Newton Law. If $\pm{\mathbf v}_0\in {\real}^3$ are the position coordinates of the centres, $\mm_\pm$ their masses;
${\mathbf v}$, with  ${\mathbf v}\ne \pm{\mathbf v}_0$, the position coordinate of the moving particle; ${\mathbf u}=\dot{\mathbf v}$ its velocity, and $1$ its mass,
the Hamiltonian of the system ({\it Euler Hamiltonian}) is
\beqa{2Cold}\JJ=\frac{\|{\mathbf u}\|^2}{2}-\frac{\mm_+}{\|{\mathbf v}+{\mathbf v}_0\|}-\frac{\mm_-}{\|{\mathbf v}-{\mathbf v}_0\|}\ ,\eeqa
with $\|\cdot\|$ being the Euclidean distance in ${\real}^3$.
Euler showed~\cite{jacobi09} that $\JJ$ exhibits  $2$ independent first integrals, in involution. One of these first integrals is the projection
\beqa{Theta}\Theta= {\mathbf M}\cdot \frac{{\mathbf v}_0}{\|{\mathbf v}_0\|}\eeqa
of the angular momentum ${\mathbf M}={\mathbf v}\times {\mathbf u}$ of the particle along the direction ${\mathbf v}_0$. It
 is not specifically due to the Newtonian potential, but, rather, to  its invariance  by rotations around the axis ${\mathbf  v}_0$. For example, it  persists if the Newtonian potential is replaced with a $\a$--homogeneous one.
The existence of the following constant of motion, which we shall refer to as {\it Euler integral}:

\beqa{G1}\EE=\|{\mathbf v}\times {\mathbf u}\|^2+({\mathbf v}_0\cdot {\mathbf u})^2+2 {\mathbf v}\cdot {\mathbf v}_0\left(\frac{\mm_+}{\|{\mathbf v}+{\mathbf v}_0\|}-\frac{\mm_-}{\|{\mathbf v}-{\mathbf v}_0\|}\right)\eeqa
is pretty specific of $\JJ$. As observed in~\cite{bekovO78}, in the limit of merging centres, i.e., ${\mathbf v}_0=\mathbf 0$, $\JJ$ reduces to the Kepler Hamiltonian, and  $\EE$ to the squared length of the angular momentum of the moving particle. {Note however that  $\JJ$ reduces to the Kepler Hamiltonian also when $\mm_-$ (or $\mm_+$) vanishes. The limiting value of $\EE$ for  this case is precisely what is studied in the present paper.}

The formula in~\equ{G1} is  not easy\footnote{See however~\cite{dullinM16} for a formula related to~\equ{G1}.} to be found in the literature.  There is a  classical argument  of separation of variables (which we shall recall in Section~\ref{The classical integration of the two--centre problem}) which, besides showing the integrability of~\equ{2Cold}, also can be used to derive~\equ{G1}. Such argument, however, does not provide a complete outline of the problem, since, as a matter of fact, leaves important questions unanswered, like, as an example,   the existence of action--angle coordinates, of periodic orbits, the complete picture of the bifurcation diagram.
Because of this, the problem has received, in the last decades, a renewed interest and noticeable papers appeared~\cite{waalkensDR06, biscaniI16, dullinM16, boscagginDT17}. A common ingredient of the mentioned literature is
a separation--like change of coordinates, possibly combined with a ``regularising'' change of time, which  allows, following Euler's ideas, to decouple the Hamiltonian.

Our approach to the problem is, in a sense, affected by methods of perturbation theory, and goes as follows.  We do not use decoupling coordinates and, {for conveniency}, begin with a situation where the  attracting centres are in a   ``asymmetric'' position. Namely, in place of~\equ{2Cold}, we write
\beqa{newH2C}
 \JJ=\frac{\|{{\mathbf y}}\|^2}{2 \mm}
-\frac{\mm\cM}{\|{\mathbf x}\|}-\frac{\mm\cM'}{\|{\mathbf x}'-{{\mathbf x}}\|}
\eeqa
Here, $\mm$ is the mass of the moving particle, $({\mathbf y}, {\mathbf x})$, with ${\mathbf y}=\mm\dot{\mathbf x}$ are its impulse--position coordinates, and $\cM$, $\cM'$ are the masses of the two attracting centres, posed at $\mathbf 0$, $\mathbf x'$, respectively.

In this case,  as we shall show below,  apart for a negligible additive term,
  its Euler integral takes the expression
\beqa{ENEW}\EE=\|{\mathbf M}\|^2-{\mathbf x}'\cdot {\mathbf L}+\mm^2{\mathcal  M}'\frac{({\mathbf x}'-{\mathbf x})\cdot {\mathbf x}'}{\|{\mathbf x}'-{\mathbf x}\|}\eeqa
where
\beqa{CL}
{\mathbf M}:={\mathbf x}\times {\mathbf y}\ ,\qquad {\mathbf L}:={{\mathbf y}}\times{\mathbf M}-\mm^2{\mathcal  M}\frac{{{\mathbf x}}}{\|{{\mathbf x}}\|}=\mm^2\cM\,\ee {\mathbf P}\eeqa
are the  {\it angular momentum} and the {\it eccentricity vector}   associated to the Kepler Hamiltonian
\beqa{kepler}\JJ_0:=\frac{\|{{\mathbf y}}\|^2}{2 \mm}
-\frac{\mm\cM}{\|{\mathbf x}\|}\eeqa
with $\ee$ and ${\mathbf P}$ being the eccentricity and the perihelion direction ($\|{\mathbf P}\|=1$).
{With these notations,} $\JJ$ reduces to a Kepler Hamiltonian  either for ${\mathbf x}'=\mathbf 0$, in which case, as in the symmetric case above, $\EE$ reduces to $\|{\mathbf M}\|^2${;} {or} for $\cM'=0$. In  {the latter} case, $\JJ$ and $\EE$ become, respectively,
$\JJ_0$ in~\equ{kepler}
and
\beqa{EEE}\EE_0&=&\|{\mathbf M}\|^2-{\mathbf x}'\cdot {\mathbf L}\eeqa
which{, as expected,} is a combination of first integrals of $\JJ_0$.

 {The second (and main)} difference with the traditional approach to the problem is that, as mentioned, we {\it do not use elliptic coordinates}. More closely to a perturbative point of view, we use a special {\it partial Kepler map} which reduces $\JJ$ to a two--degrees of freedom Hamiltonian.   We call {\it partial Kepler map}   any canonical map
\beqa{C}\cC:\qquad(\L,\ell,u,v)\in \cA\times\torus\times V\to (\underline{\mathbf y}, \underline{\mathbf x})=({\mathbf y}', {\mathbf y}, {\mathbf x}', {\mathbf x})\in (\real^3)^4\eeqa
 where $\cA$ is a domain\footnote{{By ``domain'' we mean an open and connected set in ${\mathbb K}=\real^m, \complex^m$.}} in $\real$,   $V$ is a domain in $\real^{10}$, $\torus:=\real/(2\p\integer)$ is the standard torus, $(u, v)=\big((u_1, u_2, u_3, u_4, u_5)$, $(v_1, v_2, v_3, v_4, v_5)\big)$,
 which ``preserves the standard two--form'':
\beqano d{\mathbf y}'\wedge d{\mathbf x}'+d{\mathbf y}\wedge d{\mathbf x}=d\L\wedge d\ell+d u\wedge d v\eeqano
and ``integrates the Keplerian motions of $({\mathbf y}, {\mathbf x})$'':
 \beqa{2B}
\left(\frac{\|{\mathbf y}\|^2}{2\mm}-\frac{\mm\cM}{\|{\mathbf x}\|}\right)\circ\cC=
-\frac{\mm^3\cM^2}{2\L^2}\ ,\eeqa
where $\mm$, $\cM$ are fixed ``mass parameters''. Of course, we have assumed that the image of $\cC$ in~\equ{C}
is a domain  where the left hand side of~\equ{2B} takes negative values.
We consider the Lagrange average of the Newtonian potential in~\equ{newH2C} written in terms of $\cC$, namely, the function
 \beqa{h1} \UU(\L,u,v):=-\frac{\mm\cM'}{2\p}\int_\torus \frac{d\ell}{\|{\mathbf x}'(\L,\ell, u, v)-{\mathbf x}(\L,\ell, u, v)\|}\eeqa We call such function {\it partially averaged Newtonian potential}. This function has been investigated in~\cite{pinzari19}. We recall the  main results of that analysis.

\nl
(i)  As a six degrees of freedom Hamiltonian,   $\UU$ is integrable by quadratures for possessing, besides itself, five independent and commuting first integrals
 which Poisson--commute with it. These are:
 \begin{itemize}
\item[{\rm $\rm I_1$:=}] the semi--major axis action $\L:=\mm\sqrt{\cM a}$;
\item[{\rm $\rm I_2$:=}] the Euclidean length $\rr:=\|{\mathbf x}'\|$;
\item[{\rm $\rm I_3$:=}]  the Euclidean length of  the  total angular momentum ${\mathbf C}:={\mathbf x}'\times {\mathbf y}'+{\mathbf x}\times {\mathbf y}$, with  ``$\times$'' denoting skew--product;
\item[{\rm $\rm I_4$:=}]  its third component $\ZZ:={\mathbf C}\cdot {\mathbf k}$, where $({\mathbf i}, {\mathbf j}, {\mathbf k})$ is a prefixed orthonormal frame;
\item[{\rm $\rm I_5$:=}] the projection  of  ${\mathbf M}={\mathbf x}\times {\mathbf y}$ along the direction ${\mathbf x}'$, defined as in~\equ{Theta}, with ${\mathbf v}_0$ replaced by  ${\mathbf x}'$.
\end{itemize}

\nl
(ii)  Besides  with $\II_1$, $\cdots$, $\II_5$, $\UU$ also Poisson--commutes with the function $\EE_0$ in~\equ{EEE}, which turns to be independent of, and commuting with, $\II_1$, $\cdots$, $\II_5$.

\nl
(iii)  There is a special partial Kepler map, which we denote as $\cK$,  which includes $\II_1$, $\cdots$, $\II_5$
 among its coordinates. Written in  terms of $\cK$,  $\UU$ and $\EE_0$  depend only on {\it one} coordinate couple of canonical coordinates, which we denote  $(\GG, {\rm g}$). Here,  $\GG=\|\mathbf x\times \mathbf y\|$ and ${\rm g}$ defines the direction of ${\mathbf P}$ in a suitable reference frame.

 \nl
(iv) The most remarkable property   (which ~\cite{pinzari19} has been called {\it renormalizable integrability}) is that $\UU$ depends on the coordinates $(\GG, {\rm g})$ only as a function of the function $\EE_0$ in~\equ{EEE}.

\nl
(v)  As a consequence of (iv), apart for certain particular initial data that can be described in closed form, the motions of the only coordinate couple of $\cK$ that moves  are the same, whether under $\UU$ or $\EE_0$, up to an unessential change of time. 
In other words, the phase portraits of the functions $\EE_0$ and $\UU$ expressed in terms of $\cK$ coincide.  The utility of this assertion relies on the fact that,    while $\UU$ is defined after a quadrature,  the expression of
       $\EE_0$ in terms  of $\cK$ is  very simple.

 The purpose of this paper is  to discuss the  dynamical aspects.
More precisely, this paper is organised as follows. In Sections~\ref{The classical integration of the two--centre problem} and~\ref{The ``asymmetric'' case}, we recall the classical integrability argument of the Hamiltonian~\equ{2Cold} and derive the formulae in~\equ{G1},~\equ{ENEW} and~\equ{EEE}. In Section~\ref{coordinates}, we define the $\cK$--coordinates and provide the expressions of $\JJ$, $\UU$ and $\EE$ in their terms. In particular, the function $\UU$, expressed in terms of $\cK$, is one--dimensional. In Section~\ref{The secular Euler Hamiltonian} we review the concept of the mentioned renormalizable integrability. In Section~\ref{Dynamical consequences} we discuss the dynamical consequences in the particular  case of the planar problem.  In Section~\ref{An application to the three--body problem (sketch)} we outline a possible  application of the results of the paper to the three--body problem, deferring  the complete analysis to a next paper.
 For definiteness, and by the author's tastes, we just look at  the  ``full'' problem.  The author is not aware whether addressing the same question to the ``restricted''  problem would be simpler.

\subsection{The classical integration of the two--centre problem} \label{The classical integration of the two--centre problem}

Let $\JJ$ be as in~\equ{2Cold}.
After fixing a reference frame with the third axis in the direction of ${\mathbf v}_0$ and denoting as $(v_1,  v_2,  v_3)$ the coordinates of $ {\mathbf v}$ with respect to such frame,
one introduces
the so--called  ``elliptic coordinates''
 \beqa{lambdabeta} \l=\frac{1}{2}\left(\frac{\rr_+}{\rr_0}+\frac{\rr_-}{\rr_0}\right)\ ,\quad \b=\frac{1}{2}\left(\frac{\rr_+}{\rr_0}-\frac{\rr_-}{\rr_0}\right)\ ,\quad \o:=\arg{(-v_2, v_1)}\eeqa
 where we have let, for short,
\beqano\rr_0:=\|{\mathbf v}_0\|\ ,\quad \rr_\pm:=\|{\mathbf v}\pm {\mathbf v}_0\|\ .\eeqano

Regarding $\rr_0$ as a fixed external parameter and calling $p_\l$, $p_\b$, $p_\o$ the generalized momenta associated to $\l$, $\b$ and $\o$, it turns out that
the Hamiltonian~\equ{2Cold}, written in the coordinates $(p_\l, p_\b, \l, \b)$ is independent of $\o$  and has the expression
\beqa{H2C}\JJ(p_\l, p_\b, p_\o, \l, \b, \rr_0)&=&\frac{1}{\l^2-\b^2}\Big[\frac{p^2_\l(\l^2-1)}{2 \rr_0^2}+\frac{p^2_\beta(1-\beta^2)}{2 \rr_0^2}\nonumber\\
&&+\frac{p_\o^2}{2 \rr_0^2}\big(\frac{1}{1-\beta^2}+\frac{1}{\l^2-1}\big)-\frac{(\mm_++\mm_-)\l}{\rr_0^2}\nonumber\\
&&+\frac{(\mm_+-\mm_-)\beta}{\rr_0^2}\Big]\ .
\eeqa
It follows that the solution of ``Hamilton--Jacobi'' equation
\beqa{HJNEW}\JJ(W_\l, W_\b, p_\o, \l, \b, \rr_0)=h\eeqa
can be searched of the form
\beqano W(\l, \b, p_\o, \rr_0, h)=W^\ppu(\l,  p_\o, \rr_0, h)+W^\ppd(\b, p_\o, \rr_0, h)\eeqano
as~\equ{HJNEW} separates completely as
\beqa{split}{\cF}_1(W^\ppu_\l,\l,{p_\o}, \rr_0, h)+{\cF}_2(W^\ppd_\b,\b,{p_\o}, \rr_0, h)=0\eeqa
with
\beqano
&&{\cF}_1(p_\l,\l,{p_\o}, \rr_0, h)=p^2_\l(\l^2-1)+\frac{p_\o^2}{\l^2-1}-2(\mm_++\mm_-)\l-2 \rr_0^2\l^2h\nonumber\\
&&{\cF}_2(p_\b,\b,{p_\o}, \rr_0, h)=p^2_\beta(1-\beta^2)+\frac{p_\o^2}{1-\beta^2}+2(\mm_+-\mm_-)\beta+2 \rr_0^2\beta^2h\ .
\eeqano

The identity~\equ{split} implies that
there must exist a function $\EE$, which we call {\it Euler integral}, depending on  $({p_\o}, \rr_0, h)$ only,  such that
\beqano{\cF}_\l(p_\l,\l,{p_\o}, \rr_0, h)=-{\cF}_\b(p_\l,\l,{p_\o}, \rr_0, h)=\EE({p_\o}, \rr_0, h)\quad \forall\ (p_\l, p_\b, \l, \b)\ .\eeqano It is given by
\beqa{Es}
\EE=\frac{1}{2}({\cF}_\b-{\cF}_\l)&=&\frac{p^2_\beta}{2}(1-\beta^2)-\frac{p^2_\l}{2}(\l^2-1)+\frac{p_\o^2}{2}\big(\frac{1}{1-\beta^2}-\frac{1}{\l^2-1}\big)\nonumber\\
&+&\mm_+(\l+\beta)+\mm_-(\l-\beta)+\rr_0^2(\l^2+\beta^2)h\ .
\eeqa
 We now check that the function $\EE$, written in the initial coordinates $({\mathbf u}$, ${\mathbf v})$,  coincides with~\equ{G1}.
To this end, we replace the  coordinates  $({\mathbf u}$, ${\mathbf v})$ with the ``canonical spherical coordinates relatively to ${\mathbf v}_0$'', which we denote as \beqano\cD_{\mathbf v_0}=(\Theta, {\MM}, \RR,\vartheta,  {m}, \rr)\ .\eeqano
Their definition is as follows. Given three vectors  ${\mathbf n}_1$, ${\mathbf n}_2$, ${\mathbf b}\in {\real}^3$, with  ${\mathbf n}_1$, ${\mathbf n}_2\perp$, ${\mathbf b}$, let
$\a_{{\mathbf b}}({\mathbf n}_1,{\mathbf n}_2)$ denotes the oriented angle defined by the ordered couple $({\mathbf n}_1,{\mathbf n}_2)$, relatively to the positive verse established by ${\mathbf b}$. Let
\beqano{\mathbf M}:={\mathbf v}\times {\mathbf u}\ ,\quad {\mathbf n}_0:={\mathbf v}_0\times{\mathbf M}\ ,\quad {\mathbf n}:={\mathbf M}\times{\mathbf v}\ .\eeqano
Then we define the coordinates
$\cD_{\mathbf v_0}$
via the formulae
\beqa{delaunay}
 \arr{\Theta:=\frac{{\mathbf M}\cdot {\mathbf v}_0}{\|{\mathbf v}_0\|}\\
{\MM}:=\|{\mathbf M}\|\\
\RR:=\frac{{\mathbf u}\cdot {\mathbf v}}{\|{\mathbf v}\|}}\qquad
\arr{
\dst\vartheta:=\a_{{\mathbf v}_0}({\mathbf i}, {\mathbf n}_0)\\
\dst {m}:=\a_{{\mathbf M}}({\mathbf n}_0, {\mathbf v})\\
\dst  \rr:=\|{\mathbf v}\|
}
\eeqa
As it is well known, the coordinates $\cD_{{\mathbf v}_0}$ are homogeneous--canonical\footnote{Namely, they leave the standard 1--form unvaried: \beqano{\mathbf u}\cdot d{\mathbf v}:=\sum_{i=1}^3 u_i dv_i=\Theta d\vartheta+\MM d m+\RR d\rr\ .\eeqano}.

Since $\Theta$ is a first integral to $\JJ$,  this Hamiltonian  will depend only on  the four coordinates
$({\MM}, \RR, {m}, \rr)$
while the action $\Theta$ will play the r\^ole of an ``external parameter'', together with
 $\rr_0$. Using such coordinates, $\JJ$ becomes
\beqa{hsDel}\JJ=\frac{\RR^2}{2}+\frac{{\MM}^2}{2\rr^2}-\frac{\mm_+}{\rr_+}-\frac{\mm_-}{\rr_-}\eeqa
with
\beqano{\rr_\pm:=\sqrt{\rr_0^2\mp2 \rr_0 \rr \sqrt{1-\frac{\Theta^2}{{\MM}^2}}\cos{m} +\rr^2}}\ .\eeqano
Combining this and
~\equ{lambdabeta},
one obtains
\beq{r}{\rr=\rr_0\sqrt{\l^2+\beta^2-1}\qquad {m}=\cos^{-1}\Big(-\frac{\l\beta}{\sqrt{\l^2+\beta^2-1}\sqrt{1-\frac{\Theta^2}{{\MM}^2}}}\Big)}\ .\eeq
The use of the associated generating function
\beqano
S({\MM}, {\Theta}, \l, \b)&=&\RR \rr_0\sqrt{\l^2+\beta^2-1}\nonumber\\
&+&\int^{\MM}  \cos^{-1}\left(-\frac{\l\beta}{\sqrt{\l^2+\beta^2-1}\sqrt{1-\frac{\Theta^2}{{{\MM}'}^2}}}\right)d{\MM}'
\eeqano
allows to find the generalized impulses $\ovl p_\l$, $\ovl p_\b$ associated to
$\l$, $\b$ as
\beqano\arr{
\dst \ovl p_\l=\frac{\rr_0\l \RR}{\sqrt{\l^2+\beta^2-1}}-\frac{\beta\sqrt{(1-\beta^2)(\l^2-1){\MM}^2-(\l^2+\beta^2-1)	\Theta^2}}{(\l^2+\beta^2-1)(\l^2-1)}\\
\dst \ovl p_\beta=\frac{\rr_0\beta \RR}{\sqrt{\l^2+\beta^2-1}}+\frac{\l \sqrt{(1-\beta^2)(\l^2-1){\MM}^2-(\l^2+\beta^2-1)	\Theta^2}}{(\l^2+\beta^2-1)(1-\beta^2)}
}
\eeqano
We solve for $\RR$, ${\MM}^2$:
\beqano
\arr{\dst \RR=\frac{\l(\l^2-1)\ovl p_\l+\beta(1-\beta^2)\ovl p_\beta}{\rr_0(\l^2-\beta^2)\sqrt{\l^2+\beta^2-1}}\\
\dst {\MM}^2=\frac{(\l \ovl p_\beta-\beta \ovl p_\l)^2(\l^2-1)(1-\beta^2)}{(\l^2-\beta^2)}+\frac{\l^2+\beta^2-1}{(1-\beta^2)(\l^2-1)}\Theta^2}
\eeqano
Using these formulae and the~\equ{r} inside the Hamiltonian~\equ{hsDel}, we find exactly the expression in~\equ{H2C}, with $ p_\l$, $p_\b$, $p_\o$ replaced by $\ovl p_\l$, $\ovl p_\b$, $\Theta$. Therefore,
the Euler integral will be  as in~\equ{Es}, with the same substitutions.
After some
 elementary computation, we find that the $\EE$ has, in terms of $\cD_{{\mathbf v}_0}$, the expression
\beqano
\EE&=&{\MM}^2+\rr_0^2(1-\frac{\Theta^2}{{\MM}^2})(-\RR\cos{m}+\frac{{\MM}}{\rr}\sin{m})^2\nonumber\\
&&- 2 \rr \rr_0\cos{m}\sqrt{1-\frac{\Theta^2}{{\MM}^2}}\left(\frac{\mm_+}{\rr_+}-\frac{\mm_-}{\rr_-}\right)
\eeqano
with $\rr_\pm$ as in~\equ{lambdabeta}. Turning back to the coordinates ${\mathbf u}$, ${\mathbf v}$ via~\equ{delaunay}, one sees that $\EE$ has the expression in~\equ{G1}.

\subsection{The ``asymmetric'' case}\label{The ``asymmetric'' case}
We prove that, if the  two attracting centres are  posed in ``asymmetric'' positions with respect to a prefixed reference frame, namely, the Euler Hamiltonian is written in the form~\equ{newH2C}, then  its Euler integral takes the expression in Eqs.~\equ{ENEW},~\equ{CL},~\equ{EEE}, apart for a negligible additive term.
 To this end, we let
\beqano\widehat{\JJ}(\widehat{\mathbf y},\widehat{\mathbf x}, \widehat{\mathbf x}'):=\frac{1}{\mm}{\JJ}(\mm \widehat{\mathbf y}, \widehat{\mathbf x}, \widehat{\mathbf x}')=\frac{\|  {\widehat{\mathbf y}}\|^2}{2}-\frac{{\mathcal  M}}{\|{\widehat{\mathbf x}}\|}-\frac{\cM'}{\|{\widehat{\mathbf x}}-{\widehat{\mathbf x}'}\|}\eeqano
and then we change, canonically,
\beqano\widehat{\mathbf x}'=2{\mathbf v}_0\ ,\quad \widehat{\mathbf x}={\mathbf v}_0+{\mathbf v}\ ,\quad \widehat{\mathbf y}'=
\frac{1}{2}({\mathbf u}_0-{\mathbf u})
\ ,\quad  \widehat{\mathbf y}={\mathbf u} \eeqano
(where $\widehat{\mathbf y}'$, $\widehat{\mathbf u}_0$ denote the generalized impulses conjugated to $\widehat{\mathbf x}'$, $\widehat{\mathbf v}_0$, respectively)
we reach the Hamiltonian $\JJ$ in~\equ{2Cold}, with $\mm_+={\mathcal  M}$, $\mm_-=\cM'$.  Turning back with the transformations, one sees that the function $\EE$ in~\equ{G1} takes the expression
\beqano
\frac{\EE}{\mm}&:=&\frac{1}{{\rm m}}\Big\|\Big({\mathbf x}-\frac{{\mathbf x}'}{2}\Big)\times {\mathbf y}\Big\|^2+\frac{1}{4{\rm m}}({\mathbf x}'\cdot {\mathbf y})^2\nonumber\\
&+&\mm{\mathbf x}'\cdot\Big({\mathbf x}-\frac{{\mathbf x}'}{2}\Big)\Big(\frac{{\mathcal  M}}{\|{\mathbf x}\|}-\frac{\cM' }{\|{\mathbf x}'-{\mathbf x}\|}\Big)\ .
\eeqano
After multiplying  by $\mm$, we rewrite the latter integral as
\beqano
\EE=\EE_0+\EE_1+\EE_2\eeqano
with
\beqano
&&\EE_0:=\|{\mathbf M}\|^2-{\mathbf x}'\cdot {\mathbf L}\ ,\quad
\EE_1:= \mm^2{\mathcal  M}'\frac{({\mathbf x}'-{\mathbf x})\cdot {\mathbf x}'}{\|{\mathbf x}'-{\mathbf x}\|}\nonumber\\
&& \EE_2:=\mm\frac{\|{\mathbf x}'\|^2}{2}\left(\frac{\|{{\mathbf y}}\|^2}{2 \mm}
-\frac{\mm\cM}{\|{\mathbf x}\|}-\frac{\mm\cM'}{\|{\mathbf x}'-{{\mathbf x}}\|}\right)\eeqano
where ${\mathbf M}$, ${\mathbf L}$
are as  in~\equ{CL}.
  Since $\EE_2$ is  itself an integral for ${\JJ}$, we can neglect it and rename
    \beqa{cal G new}
\EE:=\EE_0+\EE_1\eeqa
the Euler integral to ${\JJ}$.

\section{$\cK$ coordinates}\label{coordinates}

We describe a set of canonical
coordinates, which we denote as $\cK$, which we shall use for our analysis of the Euler Hamiltonian~\equ{newH2C}.

We consider, in the region of  phase space where $\JJ_0$ in~\equ{kepler} takes negative values, the ellipse with initial datum $({\mathbf y}, {\mathbf x})$. Denote as:
 \begin{itemize}
 \item[{\tiny\textbullet}]
 $a$  the {\it semi--major axis};
  \item[{\tiny\textbullet}] ${\mathbf P}$, with $\|{\mathbf P}\|=1$, the direction of perihelion, assuming the ellipse is not a circle;
    \item[{\tiny\textbullet}] $\ell$: the mean anomaly, defined, mod $2\p$, as the area of the elliptic sector spanned by ${\mathbf x}$ from ${\mathbf P}$, normalized to $2\p$.
 \item[{\tiny\textbullet}]{Given ${\mathbf u}$, ${\mathbf v}$, ${\mathbf w}\in \real^3$, with ${\mathbf u}$, ${\mathbf v}$, $\perp{\mathbf w}$, we denote as $\a_{w}({u},{v})$  the oriented angle  ${\mathbf u}$ to ${\mathbf v}$, with respect to the counterclockwise orientation established by ${\mathbf w}$.}
    \end{itemize}



We fix an arbitrary (``inertial'') frame
\beqano\FF_0:\qquad {\mathbf i}_0=\left(
\begin{array}{lll}
1\\
0\\
0
\end{array}
\right)\ ,\qquad {\mathbf j}_0=\left(
\begin{array}{lll}
0\\
1\\
0
\end{array}
\right)\ ,\qquad {\mathbf k}_0=\left(
\begin{array}{lll}
0\\
0\\
1
\end{array}
\right)\eeqano in $\real^3$,
and denote as
\beqano{\mathbf M}={\mathbf x}\times {\mathbf y}\ ,\quad {\mathbf M}'={\mathbf x}'\times {\mathbf y}'\ ,\quad {\mathbf C}={\mathbf M}'+{\mathbf M}\ ,\eeqano
where ``$\times$'' denotes skew--product in ${\real}^3$.
Observe the following relations
\beqa{orthogonality}{{\mathbf x}'}\cdot{{\mathbf C}}={{\mathbf x}'}\cdot{\big({\mathbf M}+{\mathbf M}'\big)}={{\mathbf x}'}\cdot{{\mathbf M}}\ ,\qquad {\mathbf P}\cdot {\mathbf M}=0\ ,\quad \|{\mathbf P}\|=1\ .\eeqa
Assume that the ``nodes''
\beqa{nodes}&&{\mathbf i}_1:={\mathbf k}_0\times {\mathbf C}\ ,\quad {\mathbf i}_2:={\mathbf C}\times {\mathbf x}'\ ,\quad {\mathbf i}_3:={\mathbf x}'\times {\mathbf M}\eeqa
do not vanish.
We   define the coordinates \beqano{\cK}=(\ZZ, \CC, \Theta, \GG, \RR  , \L, \zeta, g, \vartheta, {\rm g}, \rr , \ell)\eeqano
via the following formulae{:}
\beqa{belline}
\arr{
 \dst {\rm Z}:={\mathbf C}\cdot {\rm k}\\ \\
 \dst \CC:=\|{\mathbf C}\|\\ \\
  \dst  \RR:=\frac{{\mathbf y}'\cdot {\mathbf x}'}{\|{\mathbf x}'\|}\\\\
\dst  \L={\rm m}\sqrt{{\rm M} a}\\\\
 \dst\GG:=\|{\mathbf M}\|\\\\
 \dst \Theta:=\frac{{\mathbf M}\cdot {\mathbf x}'}{\|{\mathbf x}'\|}\\
 }\qquad\qquad \arr{
 \dst {\rm z}:=\a_{{\rm k}}({\mathbf i}, {\mathbf i}_1)\\ \\
 \dst g:=\a_{{\mathbf C}}({\mathbf i}_1, {\mathbf i}_2)\\ \\
  \dst \rr:=\|{\mathbf x}'\|\\ \\
\ell:=\textrm{\rm mean anomaly of } {\mathbf x}\ {\rm on\ \mathbb E}\\ \\
{\rm g}:=\a_{{\mathbf M}}({\mathbf i}_3, {\mathbf M}\times {\mathbf P})\\\\
  \dst\vartheta:=\a_{{\mathbf x}'}({\mathbf i}_2, {\mathbf i}_3)
 }
\eeqa


    The canonical character of ${\cK}$ follows from~\cite{pinzari13}. Indeed, in~\cite{pinzari13}, we considered a set of coordinates for the three--body problem\footnote{An extension to the case of an arbitrary number of planets has been successively worked out in~\cite{pinzari18}.}, that here\footnote{$(\ZZ, \CC, \RR, \ovl\RR, \Theta, \Phi, \zz, g, \rr, \ovl\rr, \vartheta, \varphi)$ are called $(\CC_3, \GG, \RR_1, \RR_2, \Theta, \Phi_2, \zeta, {\mathbb g}, \rr_1, \rr_2, \vartheta, \varphi_2)$ in~\cite{pinzari13}.} we denote  as $\cP=(\ZZ, \CC, \RR, \ovl\RR, \Theta, \Phi, \zz, g, \rr, \ovl\rr, \vartheta, \varphi)$, that are related to ${\cK}$ above via the canonical change
\beqa{planar Delaunay}\cD_{e\ell, \rm pl}:\quad (\L,\GG, \ell, {\rm g})\to (\ovl\RR, \ovl\Phi, \ovl\rr, \ovl\f)\eeqa
     usually referred to as {\it planar Delaunay map}, defined as
\beqa{p coord***}
 \arr{
 \ovl\RR=\frac{\mm^2{\mathcal  M}}{\L} \frac{\ee\sin\xi}{1-\ee\cos\xi}\\
\ovl\Phi=\GG
 }\qquad  \arr{
\ovl\rr=a(1-\ee\cos\xi)\\
\ovl\f=\n+{\rm g}-\frac{\p}{2}
 }
\eeqa
where  $\xi=\xi(\L, \GG, \ell)$, $\n=\n(\L, \GG, \ell)$ are, respectively, the eccentric and the true anomaly, defined below (see~\equ{Kepler Equation},~\equ{true anomaly}). Since the map $\cD_{e\ell, \rm pl}$ in~\equ{planar Delaunay} and the coordinates $\cP$ of~\cite{pinzari13} are canonical, so is ${\cK}$. Observe, incidentally, the unusual $\frac{\p}{2}$--shift in
\equ{p coord***}, due to the fact that, according to the definitions in~\equ{belline},  ${\rm g}$ is the longitude of ${\mathbf M}\times {\mathbf P}$ in the  plane {of ${\mathbf i}_3$, ${\mathbf j}_3$, relatively to ${\mathbf i}_3$}.

\begin{remark}\rm (i) We briefly discuss the geometrical meaning of the coordinates $\cK$, deferring to~\cite{pinzari13} or~\cite[Chapter II and Appendix E]{pinzari18} for more  details.
The definitions~\equ{belline} are based on a multiple change of reference frames. Chains of reference frames  have been firstly used  by A. Deprit, in order to extend to an arbitrary number $n$ of particles the classical ``reduction of the nodes'' discovered by Jacobi in the case $n=2$~\cite{deprit83, chierchiaPi11a, jacobi1842}.
In the case of the coordinates $\cK$, we
define three orthogonal (not necessarily orthonormal) frames  $\FF_i=({\mathbf i}_i, {\mathbf j}_i, {\mathbf k}_i)$, $i=1,\ 2,\ 3$, where ${\mathbf i}_j$ are as in~\equ{nodes}, while
\beqano
&&{\mathbf k}_1:={\mathbf C}\ ,\quad {\mathbf k}_2:={\mathbf x}'\ ,\quad {\mathbf k}_3:={\mathbf M}\  ,\quad {\mathbf j}_i:={\mathbf k}_i\times {\mathbf i}_i\eeqano
The frame $\FF_1$  is also used in~\cite{deprit83} and is often referred to as the ``invariable frame'', since it does not move under the motions of a ${\rm S}\OO(3)$--invariant Hamiltonian. $\FF_2$ and $\FF_3$ are quite specific of $\cK$.
The triples  $(\ZZ, \CC, \zz)$, $(\rr, \Theta, g)$, $(\GG, \Theta, \vartheta)$ have the meaning of ``spherical coordinates'' of  ${\mathbf C}$, $\bx$, ${\mathbf M}$ relatively to  $\FF_0$, $\FF_1$, $\FF_2$, respectively.
While the triple $(\ZZ, \CC, \zz)$ also appears\footnote{The reader should beware that, even though some groups of coordinates have been already separately used in the literature, similarly to~\cite{deprit83}, the mix $\cK$  is not  a Cartesian product of canonical coordinates (namely, the standard two--form is not  the sum of forms associated to  groups of coordinates, as it happens, for example, in the case of Delaunay coordinates).} in~\cite{deprit83}, $\Theta$ and $\vartheta$ are specific of $\cK$.
 As it can be seen from the definitions~\equ{belline}, $\Theta$ measures the convex angle between $\bx'$ and ${\mathbf C}$ (or ${\mathbf M}$, by the first identity in~\equ{orthogonality}) and vanishes when the two vectors are orthogonal (in the planar case).  The angle $\vartheta$ measures the rotation  of ${\mathbf M}$ with respect to $\bx$.
    The quadruplet $(\L, \GG, \ell, {\rm g})$ is the Delaunay set of coordinates associated to $(\by, \bx)$ in  $\FF_3$. The coordinate $\RR$ measures the radial velocity of $\bx'$ (as well known, the normal velocity is measured by $\frac{\GG}{\rr}$).

        (ii) Contrary to the Jacobi reduction of the nodes, $\cK$ are well defined also in the planar case.
    In such case, they reduce to $\cK_{\rm pl}=(\CC, \RR, \L, \GG, g, \rr, \ell, {\rm g})$, and correspond to take
    the Delaunay  coordinates
     $(\L, \GG, \ell, \ovl{\rm g})$  for $(\by, \bx)$,
     the ``symplectic polar coordinates'' $(\RR, \Phi', \rr, \varphi')$ (also used in~\cite{deprit83}) for  $(\by', \bx')$ and next reduce the angular momentum via the relations
     $\Phi'=\CC-\GG$, $\f'=g$, $\ovl{\rm g}=\f'-\n(\L, \GG, \ell)+{\rm g}$.
\end{remark}

 \subsection{Expression of $\JJ$ and $\EE$ in terms of $\cK$}\label{JandE}
 Using the formulae in the previous section, we provide the expressions of $\JJ$ in~\equ{newH2C} and and $\EE$ in~\equ{EEE} in terms of ${\cK}$:
\beqa{EEJJ}
\JJ(\L, \GG, \Theta, \rr , \ell, {\rm g})&=&-\frac{\mm^3{\mathcal M}^2}{2\L^2}-\frac{\mm{\mathcal M}'}{\sqrt{{\rr }^2+2\rr  a\sqrt{1-\frac{\Theta^2}{\GG^2}} {\rm p}+{a}^2\varrho^2}}\nonumber
\eeqa\beqa{EEJJ}
&=:&\JJ_{0}+\JJ_{1}\nonumber\\
\EE(\L, \GG, \Theta, \rr , \ell, {\rm g})&=&\GG^2+\mm^2{\mathcal M}'\rr \sqrt{1-\frac{\Theta^2}{\GG^2}}\sqrt{1-\frac{\GG^2}{\L^2}}\cos{\rm g}\nonumber\\
 & +&\mm^2{\mathcal M}'\rr \frac{{\rr }+a\sqrt{1-\frac{\Theta^2}{\GG^2}} {\rm p}}{\sqrt{{\rr }^2+2\rr  a\sqrt{1-\frac{\Theta^2}{\GG^2}}{\rm p}+{a}^2\varrho^2}}\nonumber\\
&=:&\EE_{0}+\EE_{1}
\eeqa
and, if $\xi=\xi(\L, \GG, \ell)$ is the {\it eccentric anomaly}, defined as the solution of
{\it Kepler equation}
 \beqa{Kepler Equation}\xi-\ee(\L,\GG)\sin\xi=\ell\eeqa
and $a=a(\L)$ the {\it semi--major axis}; $\ee=\ee(\L, \GG)$, the {\it eccentricity} of the ellipse,
$\varrho=\varrho(\L, \GG, \ell)$, ${\rm p}={\rm p}(\L, \GG, \ell, {\rm g})$ are defined as
\beqa{p}
a(\L)&=&\frac{\L^2}{\mm^2\cM}\nonumber\\
\ee(\L, \GG)&:=&\sqrt{1-\frac{\GG^2}{\L^2}}\nonumber\\
 \varrho(\L, \GG, \ell)&:=&1-\ee(\L, \GG)\cos\xi(\L, \GG, \ell)\nonumber\\
{\rm p}(\L, \GG, \ell, {\rm g})&:=&(\cos\xi(\L, \GG, \ell)-\ee(\L, \GG))\cos{\rm g}-\frac{\GG}{\L}\sin\xi(\L, \GG, \ell)\sin{\rm g}\ .\eeqa
The angle
\beqa{true anomaly}
\n(\L, \GG, \ell):=\arg\left(\cos\xi(\L, \GG, \ell)-\ee(\L, \GG),\frac{\GG}{\L}\sin\xi(\L, \GG, \ell) \right)
\eeqa
is usually referred to as {\it true anomaly}, so one recognises that ${\rm p}(\L, \GG, \ell, {\rm g})=\varrho \cos(\n+{\rm g})$.
Observe that $\EE$ and $\JJ$ are both independent of $\CC$, $\ZZ$, $\zeta$, $\g$, $\RR$, $\vartheta$, because the conjugated coordinates to these variables are first integrals of the motion. Remark, at this respect, that: (i) the couples $(\ZZ$, $\zeta)$ and $(\CC$, $\g)$ are, simultaneously, couples of first integrals to $\JJ$ and $\EE$ and  couples of cyclic coordinates, determined by the conservation of ${\mathbf C}$ and ${\mathbf x}'$; (ii) the independence on $\vartheta$ corresponds to the invariance of $\JJ$ and $\EE$ under the one--parameter group of rotations around ${\mathbf x}'$.

The details on the derivation of the formulae in~\equ{EEJJ} may be found in Appendix~\ref{JEK}.

\section{Renormalizable integrability}\label{The secular Euler Hamiltonian}

In this section we review the property of {\it renormalizable integrability} pointed out in~\cite{pinzari19}.


 \begin{definition}\label{def: renorm integr}\rm Let $h$, $g$ be two  functions
of the form
\beqa{HJ}h(p, q, y, x)=\widehat h({\rm I}(p,q), y, x)\ ,\qquad g(p, q, y, x)=\widehat g({\rm I}(p,q), y, x)\eeqa
where
\beqa{D}(p, q, y, x)\in \cD:=\cB\times U\eeqa
with $ U\subset \real^2$, $\cB\subset\real^{2n}$ open and connected, $(p,q)=$ $(p_1$, $\cdots$, $p_n$, $q_1$, $\cdots$, $q_n)$   {conjugate} coordinates with respect to the two--form $\o=dy\wedge dx+\sum_{i=1}^{n}dp_i\wedge dq_i$ and ${\rm I}(p,q)=({\rm I}_1(p,q), \cdots, {\rm I}_n(p,q))$, with
\beqano{\rm I}_i:\ \cB\to \real\ ,\qquad i=1,\cdots, n\eeqano
pairwise Poisson commuting:
\beqa{comm}\big\{{\rm I}_i,{\rm I}_j\big\}=0\qquad \forall\ 1\le i<j\le n\qquad i=1,\cdots, n\ .\eeqa
 We say that $h$ is {\it renormalizably integrable via $g$} if there exists a function  \beqano\widetilde h:\qquad {\rm I}(\cB)\times g(U)\to \real\ , \eeqano such that
\beqa{renorm}h(p,q,y,x)=\widetilde h({\rm I}(p,q), \widehat g({\rm I}(p,q),y,x))\eeqa
for all $(p, q, y, x)\in \cD$.
\end{definition}
\begin{proposition}\label{rem}
If $h$  is renormalizably integrable via $g$, then: (i)
${\rm I}_1$, $\cdots$, ${\rm I}_n$  are first integrals to $h$ and $g$;
 (ii)
$h$ and $g$ Poisson commute. \end{proposition}

Observe that,  if $h$ is renormalizably integrable via $g$, then, generically, their respective time laws for the coordinates $(y, x)$ are the same, up to rescale the time:

\begin{proposition}\label{prop: fixed points}
{L}et $h$ be renormalizably integrable via $g$. Fix a value ${\rm I}_0$ for the integrals ${\rm I}$ and look at the motion of
 $(y, x)$ under $h$ and $g$, on the manifold ${\rm I}={\rm I}_0$. For any fixed initial datum $(y_0, x_0)$, let
 $g_0:=g({\rm I}_0, y_0, x_0)$. If
 $\o({\rm I}_0, g_0):=\partial_{g}\tilde h({\rm I}, g)|_{({\rm I}_0, g_0)}\ne 0$, the motion $(y^h(t), x^h(t))$ with initial datum$(y_0, x_0)$ under $h$ is related to the corresponding motion $(y^g(t), x^g(t))$  under  $g$ via
 \beqano y^h(t)=y^g(\o({\rm I}_0, g_0) t)\ ,\quad x^h(t)=x^g(\o({\rm I}_0,  g_0) t)\eeqano
In particular, under this condition, all the fixed points of $g$
in the plane $(y, x)$ are fixed point to $h$.
Values of $({\rm I}_0, g_0)$ for which $\o({\rm I}_0, g_0)= 0$ provide, in the plane $(y, x)$, curves of fixed points for $h$ (which are not necessarily curves of fixed points to $g$).
\end{proposition}

{We consider the $\ell$--average of the  of the function $\JJ_1$ in~\equ{EEJJ}:
 \beqa{h1} \UU(\rr, \L, \Theta, \GG, {\rm g}):=-\frac{\mm\cM' }{2\p}\int_0^{2\p}
 \frac{d\ell}{\sqrt{{\rr }^2+2\rr  a\sqrt{1-\frac{\Theta^2}{\GG^2}} {\rm p}+{a}^2\varrho^2}}
\eeqa}

We observe that $\UU$ and $\EE_0$
 have the form in~\equ{HJ}, with $\II=(\II_1, \II_2, \II_3)=(\rr, \L, \Theta)$ verifying~\equ{comm} and $(y, x)=(\GG, {\rm g})$.

\begin{proposition}\label{main prop}
$\UU$ is renormalizably integrable via ${\rm E}_0$. Namely,
there exists a function $\FF$ such that
\beqano\UU({\rm r}, \L, \Theta, {\rm G}, {\rm g})=\FF\big({\rm r}, \L, \Theta, {\rm E}_0({\rm r}, \L, \Theta, {\rm G}, {\rm g})\big)\ .\eeqano
A\footnote{We remark that $\FF$ may have several expressions, as well as $\UU$, which is defined via a quadrature.} function $\FF$ is given by
\beqano
\FF({\rm r}, \L, \Theta, {\rm E}_0)=\widetilde{\rm F}({\rm r}, a(\L), {\mathcal E}(\L,{{\rm E}_0}), {\cI}(\L,\Theta, {{\rm E}_0}))
\eeqano
where
\beqa{U} \widetilde{\rm F}(\rr, a,  \cE ,  \cI )=\frac{1}{2\p}\int_{\torus}\frac{(1- \cE \cos w)dw}{\sqrt{\rr^2+a^2-2a(\rr \cI \sin w+a \cE \cos w)+a^2 \cE ^2\cos^2w}}\  ; \eeqa
\beqa{EI} {\mathcal E}(\L,{{\rm E}_0})=\frac{\sqrt{\L^2-{{\rm E}_0}}}{\L}\qquad {\cI}(\L,\Theta, {{\rm E}_0})=\frac{\sqrt{{{\rm E}_0}-{\Theta^2}}}{\L_2}\ .
\eeqa
\end{proposition}

The results of the present section may be summarised as follows. Proposition~\ref{main prop} implies
\beqa{commutation}\Big\{\UU,\ \EE_0\Big\}=0\ ,\eeqa
but, actually, combining it with Proposition~\ref{prop: fixed points}, we have that much more is true:

\begin{itemize}
\item[(i)] If $\FF_{\EE_0}\ne 0$, the  time laws of $({\rm G}, {\rm g})$ under $\UU$ or $\EE_0$ are basically (i.e., up to a change of time) the same;

\item[(i)] Motions of $\EE_0$ corresponding to level sets for which  $\FF_{\EE_0}= 0$ are fixed points curves to $\UU$ (``frozen orbits''). In~\cite{pinzari19} we provided an example of frozen orbit of $\UU$ in the case $\d:=\frac{\rr}{a}\ll 1$;

\item[(iii)] $\UU$ and $\EE_0$ have the same action--angle coordinates.

\end{itemize}

\nl
In the next section, we investigate the dynamical properties of $\EE_0$ for the planar case ($\Theta=0$). The study of the spatial case appears  less explicit, so it is deferred to a next publication.


 \section{Dynamical properties of $\EE_0$ in the case $\Theta=0$}\label{Dynamical consequences}
In this section we focus on the dynamical consequences of renormalizable integrability.


\subsection{Phase portrait}\label{phase portrait}

We study the phase portrait of $\EE_0$ in~\equ{EEJJ} setting $\Theta=0$. In this case, we have to study the curves
\beqa{level sets}\dst\EE_0(\L, \GG, {\rm g};\rr )=\GG^2+\mm^2{\mathcal M}'\rr \sqrt{1-\frac{\GG^2}{\L^2}}\cos{\rm g}=\cE
\eeqa
in the plane $({\rm g}, \GG)$.
 To simplify notations, we
divide  this equation  by $\L^2$, and we rewrite it as
\beq{G0pl}\widehat\EE_0({\rm g}, \widehat\GG)=\widehat\GG^2+\d \sqrt{1-\widehat\GG^2}\cos{\rm g}=\widehat\cE\ , \eeq
where
\beqa{delta}
\widehat\cE:=\frac{\cE}{\L^2}\ ,\quad \widehat\GG:=\frac{\GG}{\L}\ ,\quad  \d:=\mm^2{\mathcal M}'\frac{\rr }{\L^2}=\frac{\rr}{a}\eeqa
and we  study the rescaled level sets~\equ{G0pl} in the plane $({\rm g}, \widehat\GG)$. Observe, incidentally, that the level sets
\equ{G0pl} extend also for $\d<0$, due to the symmetry of $\EE_0$ for
\beqano(\d, {\rm g})\to (-\d, \p-{\rm g})\ .\eeqano
We limit to study them in the case
 $\d> 0$.

 For  $\d\in (0,2)$, the function $\widehat\EE_0({\rm g}, \widehat\GG)$ has a minimum, a saddle and a maximum, respectively at
\beqano\hat{\mathbf P}_{-}=(\p,0)\ ,\qquad \hat{\mathbf P}_0=(0,0)\ ,\qquad  \hat{\mathbf P}_+=\left(0,\sqrt{1-\frac{\d^2}{4}}\right)\eeqano
where it takes the values, respectively,
\beqano\widehat\cE_-=-\d\ ,\qquad\widehat\cE_0=\d\ ,\quad \widehat\cE_+=1+\frac{\d^2}{4}\ .\eeqano
Thus,  the level sets in~\equ{G0pl} are non--empty only for \beqa{cal J}\widehat\cE\in 
\left[-\d, 1+\frac{\d^2}{4}\right]\ .\eeqa
We
denote as $\cS_0$,  the level set through  the saddle ${\mathbf P}_0$. When $\widehat\GG=1$, $\widehat\EE_0$ takes the value $1$ for all ${\rm g}$ and we denote as $\cS_1$ the level curve with $\widehat\cE=1$. The equations of $\cS_0$, $\cS_1$ are, respectively:
\beqa{separatrices}
&&\cS_0(\d)=\Big\{
({\rm g}, \widehat\GG):\ \widehat\GG^2+\d \sqrt{1-\widehat\GG^2}\cos{\rm g}=\d
\Big\}
\nonumber\\
&&\cS_1(\d)=\Big\{\widehat\GG=\pm 1\Big\}\cup\Big\{
\widehat\GG_1=\pm\sqrt{1-\d^2\cos^2{\rm g}}
\Big\}
\eeqa
$\cS_1$ is composed of two branches,
which will be referred to as ``horizontal'', ``vertical'', respectively,
 glue smoothly at
$(\pm \frac{\p}{2}, 1)$, with ${\rm g}$ mod $2\p$.
Note that, when $0\le \d\le 1$, the vertical branch is defined for all ${\rm g}\in\torus$; when $\d>1$,  its domain in ${\rm g}$ is made of two disjoint neighborhoods of $\pm \frac{\p}{2}$.

\begin{figure}
\includegraphics[width=0.30\textwidth]{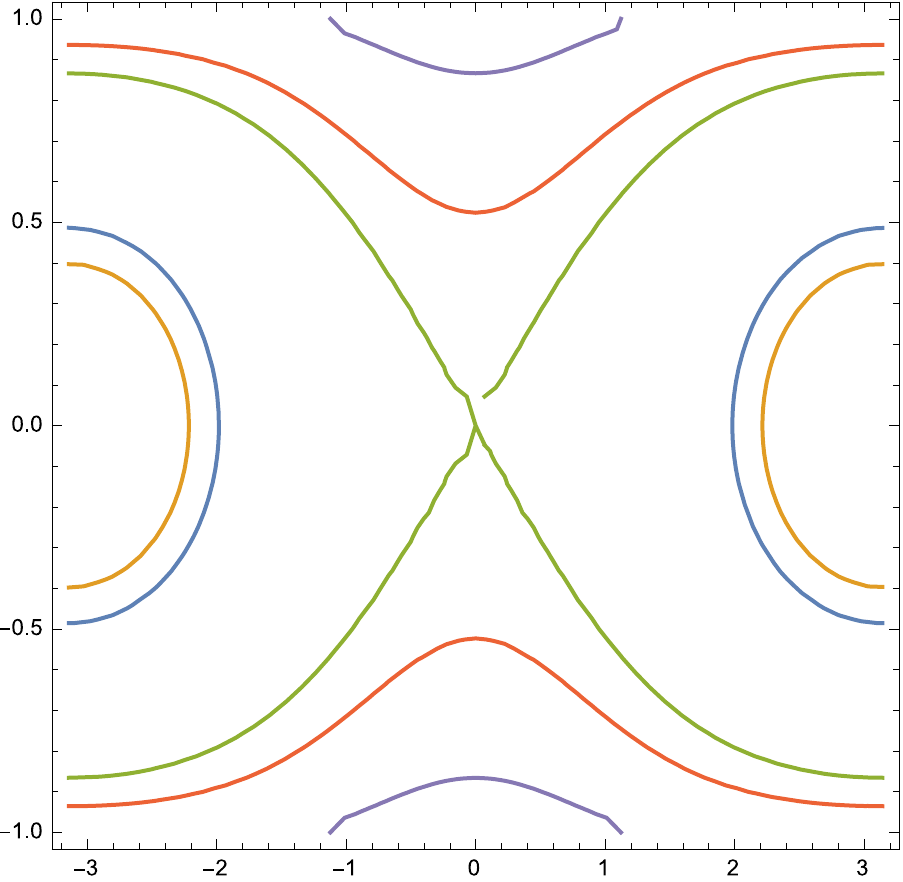}
\hfill
\includegraphics[width=0.30\textwidth]{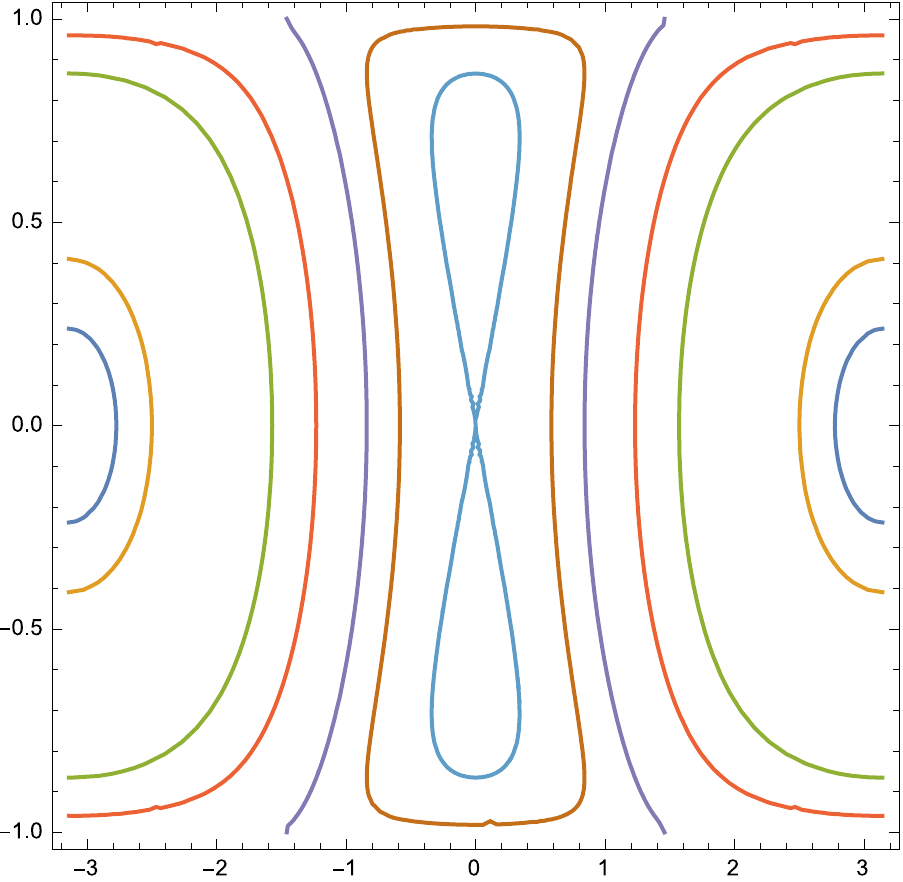}
\hfill
\includegraphics[width=0.30\textwidth]{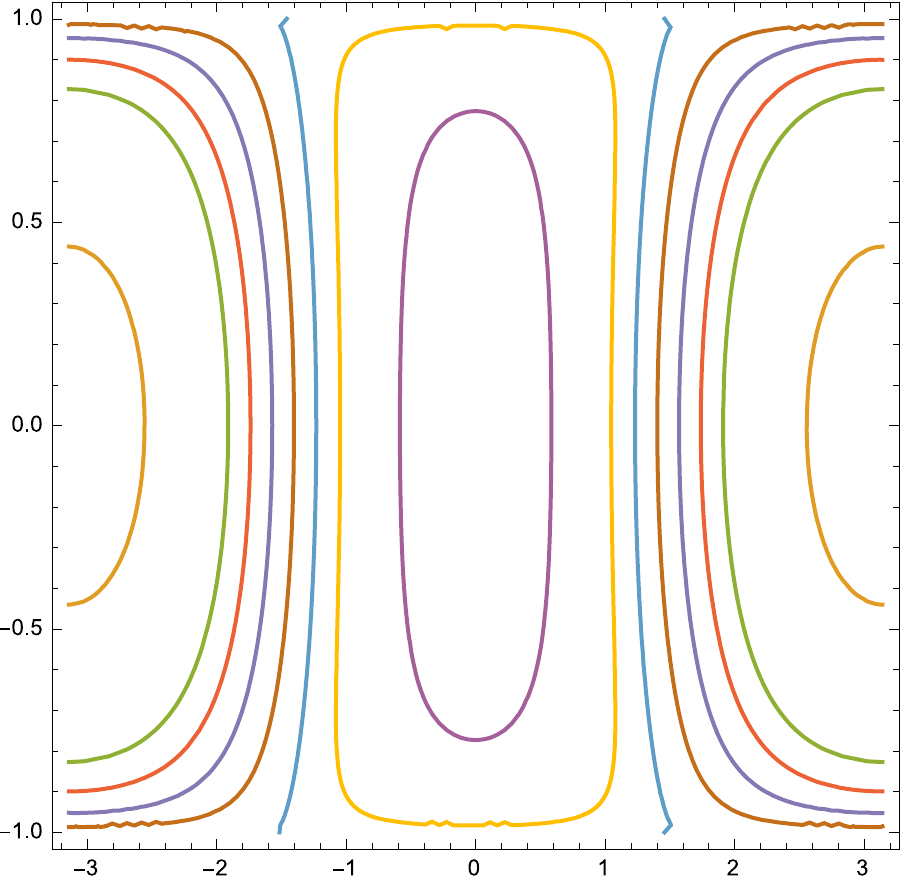}
\caption{The phase portrait of $\EE_0$ in the plane $({\rm g}, \GG)$. Left: $0<\d<1$;  Center: $0<\d<1$; Right: $\d>2$.}
\label{figure1}
\end{figure}

When $\d>2$, the saddle $\hat{\mathbf P}_0$ and its manifold $\cS_0$ do not exist, $\hat{\mathbf P}_-=(\p, 0)$ is still a minimum, while $\hat{\mathbf P}_+=(0, 0)$ becomes a maximum. The manifold $\cS_1$ still exists, with the vertical branch closer and closer, as $\d\to+\infty$, to the portion of straight ${\rm g}=\pm\frac{\p}{2}$ in the strip $-1\le \widehat\GG\le 1$.
In this case the admissible values for $\widehat\cE$ are
\beqano\widehat\cE\in 
\left[-\d, \d\right]\ .\eeqano
We now turn to the phase portrait induced by the level curves~\equ{G0pl} in the plane $({\rm g}, \GG)$.  According to the value of $\d$,  the scenario changes, as detailed below. See also Figure~\ref{figure1}.

\paragraph{Phase portrait}

\begin{itemize}\label{PP}\rm
\item[$1)$]  $0<\d\le1$
\begin{itemize}
\item[$1_1)$]  $-\d\le\widehat\cE<\d$:  the level curves encircle the minimum ${\mathbf P}_-$; 
\item[$1_2)$] $\widehat\cE=\d$ is the manifold $\cS_0$;
\item[$1_3)$]  $\d<\widehat\cE<1$:  the level curves are defined for all ${\rm g}\in \torus$;
\item[$1_4)$] $\widehat\cE=1$ is the manifold $\cS_0$;
\item[$1_5)$]   $1<\widehat\cE\le1+\frac{\d^2}{4}$:  the level curves encircle the maximum ${\mathbf P}_+$;
    \end{itemize}
    Here,  with $\d=1$  the items $1_2)$, $1_3)$ and $1_4)$, and hence the two manifolds $\cS_0$ and $\cS_1$,merge.
   With $\d=0$, $1_1)$ merges with $1_2)$; $1_5)$ merges with $1_4)$. $\cS_0$ is contracted to the ${\rm g}$--axis;  $\cS_1$ is contracted to the axis $\widehat\GG=1$; the level curves are straight lines, parallel to the ${\rm g}$--axis. The level sets reduce to $\GG=\const$.
\item[$2)$] $1<\d\le 2$
\begin{itemize}
\item[$2_1)$]  $-\d\le\widehat\cE<1$:  the level curves encircle the minimum ${\mathbf P}_-$;
\item[$2_2)$]  $\widehat\cE=1$ is the manifold $\cS_1$:
\item[$2_3)$]  $1<\widehat\cE<\d$:  the level curves encircle the saddle ${\mathbf P}_0$;
\item[$2_4)$] $\widehat\cE=\d$ is the manifold $\cS_0$;
\item[$2_5)$]  $\d<\widehat\cE\le1+\frac{\d^2}{4}$;
the level curves encircle the maximum ${\mathbf P}_+$.
    \end{itemize}
    Here, when $\d=2$,  items $2_4)$ and $2_5)$ merge and the manifold $\cS_0$ with the librations inside are contracted to the point ${\mathbf P}_0$.
\item[$3)$] $\d>2$
    \begin{itemize}
    \item[$3_1)$]  $-\d\le\widehat\cE<1$:  the level curves encircle the minimum ${\mathbf P}_-$;
\item[$3_2)$]  $\widehat\cE=1$ corresponds to the manifold $\cS_1$;
\item[$3_3)$]  $1<\widehat\cE<\d$:  the level curves encircle the maximum ${\mathbf P}_0$.
        \end{itemize}
            \end{itemize}

\begin{remark}[{\it perihelion librations}]\rm
Figure~\ref{figure1} shows that librational motions of $({\rm g}, \GG)$ about ${\mathbf P}_-=(\p, 0)$ are stable for any value of $\d$. Physically, such librations correspond to small  periodic oscillations of the perihelion of the instantaneous ellipse of $(\by,\bx)$
along\footnote{Recall that, in the planar case, the perihelion anomaly is ${\rm g}-\frac{\p}{2}$, by~\equ{belline}, relatively to ${\mathbf i}_3$. Therefore,  ${\rm g}=\p$ corresponds to ${\mathbf P}\parallel {\mathbf j}_3$.} the direction ${\mathbf j}_3$, while the ellipse squeezes to a segment before reversing its direction and  again decreasing its eccentricity. Other kind of librations of course occur in the three cases, as Figure~\ref{figure1} shows. It is  quite astonishing to see that rotational motions of the perihelion do exist only when $0<\d<1$ and only between $\cS_0$ and $\cS_1$. The singularity of this fact is even more evident if one recalls that, for arbitrarily large values of $\d$, in the planetary problem, rotational motions of the perihelia occupy a positive measure set in phase space~\cite{arnold63, laskarR95, fejoz04, chierchiaPi11b}.
\end{remark}

The computations leading to the phase portrait above are elementary.  We report here the complete discussion for $\d\in (0, 2)$. The case $\d>2$ is similar.

 Solving equation~\equ{G0pl}  for  ${\rm g}$,  we find two branches
\beqa{cosg}{\rm g}={\rm g}_\pm=\pm\cos^{-1}\left(\frac{\widehat\cE-\widehat\GG^2}{\d \sqrt{1-\widehat\GG^2}}\right)\quad \mod\ 2\p\ .\eeqa
Using
\beqa{Gpm}1-\left(\frac{\widehat\cE-\widehat\GG^2}{\d \sqrt{1-\widehat\GG^2}}\right)^2=\frac{\d^2-\widehat\cE^2-2(\frac{\d^2}{2}-\widehat\cE)\widehat\GG^2-\widehat\GG^4}{\d^2(1-\widehat\GG^2)}=\frac{(\widehat\GG^2-\widehat\GG_{-}^2)(\widehat\GG_{+}^2-\widehat\GG^2)}{\d^2(1-\widehat\GG^2)}\eeqa
with
\beqa{eq: Gpm}\widehat\GG^2_\pm&=&\widehat\cE-\frac{\d^2}{2}\pm\sqrt{\big(\widehat\cE-\frac{\d^2}{2}\big)^2
+\d^2-\widehat\cE^2
}\nonumber\\
&=&\widehat\cE-\frac{\d^2}{2}\pm\d\sqrt{1+\frac{\d^2}{4}-\widehat\cE}
\eeqa
one sees the  equality~\equ{cosg} is well defined for
 \beqa{limits}\widehat\GG_{\rm min}\le \widehat\GG\le \widehat\GG_{\rm max}\eeqa
where
\beqano\widehat\GG^2_{\rm min}:=\max\{\widehat\GG^2_{-}, 0\}\ ,\qquad \widehat\GG^2_{\rm max}:=\min\{\widehat\GG^2_{+}, 1\}\ .\eeqano
Note  that, when $\widehat\cE$ takes its maximum value $1+\frac{\d^2}{4}$, one has $\widehat\GG_+^2=\widehat\GG_-^2=1-\frac{\d^2}{4}$. Therefore,   by~\equ{cosg} and~\equ{limits}, the level set
with  $\widehat\cE=1+\frac{\d^2}{4}$
 reduces to the maximum point $(0, \pm\sqrt{1-\frac{\d^2}{4}})$.
Writing \beqano\widehat\GG_-^2=\frac{\widehat\cE^2-\d^2}{\widehat\cE-\frac{\d^2}{2}+\d\sqrt{1+\frac{\d^2}{4}-\widehat\cE}}\eeqano
and noticing that
\beqano1-\widehat\GG_+^2=1-\left(\widehat\cE-\frac{\d^2}{2}+\d\sqrt{1+\frac{\d^2}{4}-\widehat\cE}\right)=\left(\frac{\d}{2}-\sqrt{1+\frac{\d^2}{4}-\widehat\cE}\right)^2\ge 0\ ,\eeqano
one finds that
\beqa{obser}
\widehat\GG_{\rm min}=\arr{
0\qquad {\rm if}\quad -\d\le \widehat\cE\le \d\\
\widehat\GG_-\quad {\rm if}\quad \widehat\cE>\d
}\ ,\qquad \widehat\GG_{\rm max}=\widehat\GG_+\ .
\eeqa
Observe that
\beqano\lim_{\widehat\cE\to \d}\GG^2_{\rm min}=\lim_{\widehat\cE\to \d}\GG^2_{-}=0\ ,\quad \lim_{\widehat\cE\to \d}\GG^2_{\rm max}=\lim_{\widehat\cE\to \d}\GG^2_{\rm +}=\d(2-\d)\eeqano
and
\beqano\lim_{\widehat\cE\to 1}\GG^2_{\rm max}=\lim_{\widehat\cE\to 1}\GG^2_{+}=1\ ,\quad \lim_{\widehat\cE\to 1}\widehat\GG_-^2=1-\d^2\ ,\quad \lim_{\widehat\cE\to 1}\widehat\GG_{\rm min}^2=\max\{1-\d^2, 0\}\ ,\eeqano
which are obtained using
\beqano\lim_{\widehat\cE\to \d}\widehat\GG_\pm^2=\d-\frac{\d^2}{2}\pm\d\left(1-\frac{\d}{2}\right)\ ,\quad \lim_{\widehat\cE\to 1}\widehat\GG^2_\pm=1-\frac{\d^2}{2}\pm \frac{\d^2}{2}\eeqano
in turn implied  by~\equ{eq: Gpm}.

In particular, $\GG_{\rm min}$, is  continuous for $\widehat\cE=\d$.
%
%
%
%
The inequality~\equ{limits} defines a symmetric domain of $\widehat\GG$ with respect to the origin, consisting of the union
\beqa{union}\widehat\cD=\widehat\cD_-\cup\widehat\cD_+\eeqa
 of two symmetric intervals
\beqano\widehat\cD_-=\left[-\widehat\GG_{\rm max}\, -\widehat\GG_{\rm min}\right]\ ,\qquad \widehat\cD_+=\left[\widehat\GG_{\rm min}\, \widehat\GG_{\rm max}\right]\ .\eeqano
Observe that,  for  $\widehat\cE>\d$ and $\widehat\cE\ne 1$, the union~\equ{union} is disjoint, since, in this case,
$\widehat\GG_{\rm min}>0$ (se~\equ{obser}).
The functions
  ${\rm g}_\pm$ in~\equ{cosg} are even functions of $\widehat\GG$ on such symmetric domain.

We  now study the curves in~\equ{G0pl}
in the plane $({\rm g}, \widehat\GG)$, for $\widehat\cE$ as in~\equ{cal J}. By symmetry, we limit to study the behavior of  ${\rm g}_+$  for $\widehat\GG\in \widehat\cD_+$. We denote as
\beqa{ovlundlg1}\underline{\rm g}:=\cos^{-1}\left(\frac{\widehat\cE-\widehat\GG_{\rm min}^2}{\d \sqrt{1-\widehat\GG_{\rm min}^2}}\right)\ ,\quad \overline{\rm g}:=\cos^{-1}\left(\frac{\widehat\cE-\widehat\GG_{\rm max}^2}{\d \sqrt{1-\widehat\GG_{\rm max}^2}}\right)\eeqa
 the values that ${\rm g}_+$ takes at the extrema of $\widehat\cD_+$.
 The explicit value of $\underline{\rm g}$, $\ovl{\rm g}$ is
\beqa{ovlundlg2}
\underline{\rm g}=\arr{0\qquad \qquad\ \ {\rm if}\qquad\qquad \widehat\cE>\d\\
\cos^{-1}\frac{\widehat\cE}{\d}\qquad {\rm if}\quad -\d\le \widehat\cE\le \d
}\ ,\qquad \overline{\rm g}=\arr{
\p\quad {\rm if}\quad\widehat\cE<1\\
\frac{\p}{2}\quad{\rm if}\quad \widehat\cE=1\\
0\quad {\rm if}\quad\widehat\cE>1
}
\eeqa
This follows from  the definitions in~\equ{eq: Gpm} and~\equ{obser}.
In particular, from~\equ{eq: Gpm} one finds, for $(\s,\widehat\cE)\ne(+,1)$
  \beqano
\frac{\widehat\cE-\widehat\GG_{\s}^2}{\d \sqrt{1-\widehat\GG_{\s}^2}}&=&\frac{\widehat\cE-\left(\widehat\cE-\frac{\d^2}{2}+\s\d\sqrt{1+\frac{\d^2}{4}-\widehat\cE}\right)}{\d\sqrt{1-\left(\widehat\cE-\frac{\d^2}{2}+\s\d\sqrt{1+\frac{\d^2}{4}-\widehat\cE}\right)}}\nonumber\\
&=&\sign\left(\frac{\d}{2}-\s\sqrt{1+\frac{\d^2}{4}-\widehat\cE}\right)\nonumber\\
&=&\arr{+1\quad {\rm for}\quad \s=-\\
-1\quad {\rm for}\quad \s=+\ \&\ \widehat\cE<1\\
+1\quad {\rm for}\quad \s=+\ \&\ \widehat\cE>1
}
\eeqano
while, for $(\s,\widehat\cE)=(+,1)$,
 \beqano
\frac{\widehat\cE-\widehat\GG_{+}^2}{\d \sqrt{1-\widehat\GG_{+}^2}}&=&
\frac{\sqrt{1-\widehat\GG_{+}^2}}{\d}=\frac{\sqrt{1-\left(1-\frac{\d^2}{2}+\d\sqrt{1+\frac{\d^2}{4}-1}\right)}}{\d}=0\eeqano
%
%
Let us study the graph of ${\rm g}_+$ as a function of $\widehat\GG$, for $\widehat\GG\in \widehat\cD_+$.
From the formula
\beqa{g derivative}\partial_{\widehat\GG}{\rm g}_+=\frac{\widehat\GG}{\sqrt{(\widehat\GG^2-\widehat\GG^2_{-})(\widehat\GG^2_{\rm max}-\widehat\GG^2)}}\frac{2-\widehat\cE-
\widehat\GG^2}{1-\widehat\GG^2}\ .
\eeqa
%
one sees that
 $\widehat\GG=\widehat\GG_0:=\sqrt{2-\widehat\cE}\notin \widehat\cD_+$ is an extremal point, as soon as $\widehat\GG_0\in \widehat\cD_+$.
Using
\beqano
\widehat\GG_0^2-\widehat\GG_{\rm max}^2&=&2-\widehat\cE-\left(\widehat\cE-\frac{\d^2}{2}+\d\sqrt{1+\frac{\d^2}{4}-\widehat\cE}\right)\nonumber\\
&=&\sqrt{1+\frac{\d^2}{4}-\widehat\cE}\left(2\sqrt{1+\frac{\d^2}{4}-\widehat\cE}-\d\right)\nonumber\\
&=&2\frac{\sqrt{1+\frac{\d^2}{4}-\widehat\cE}}{\sqrt{1+\frac{\d^2}{4}-\widehat\cE}+\d}(1-\widehat\cE)
\eeqano
and
\beqano
\widehat\GG_0^2-\widehat\GG_{\rm min}^2&\ge&2-\widehat\cE-\left(\widehat\cE-\frac{\d^2}{2}-\d\sqrt{1+\frac{\d^2}{4}-\widehat\cE}\right)\nonumber\\
&=&\sqrt{1+\frac{\d^2}{4}-\widehat\cE}\left(2\sqrt{1+\frac{\d^2}{4}-\widehat\cE}+\d\right)\ge0\ .\eeqano
we see that
\beqano{\rm g}_0\arr{\ge \widehat\GG_{\rm max}\quad {\rm for}\quad \widehat\cE< 1\\
\in \widehat\cD_+\quad {\rm for}\quad \widehat\cE\ge 1
}\ .\eeqano
As a consequence,
\begin{itemize}
\item[{\tiny\textbullet}] For $\widehat\cE<1$, $\widehat\GG_0>\widehat\GG_{\rm max}$ and hence ${\rm g}_+$ increases,
 in $\widehat\cD_+$,
 from $\underline{\rm g}$ to $\overline{\rm g}$.
\item[{\tiny\textbullet}]   For $\widehat\cE>1$,  ${\rm g}_+$ increases from $\underline{\rm g}$  to ${\rm g}_0$ for $\widehat\GG_{\rm min}\le \widehat\GG\le \widehat\GG_0$ and decreases from ${\rm g}_0$ to  $\overline{\rm g}$, for $\widehat\GG_0\le \widehat\GG\le \widehat\GG_{\rm max}$.
\end{itemize}

\subsection{The collisional manifold and its motions}\label{The collisional manifold and its motions}

The $\EE_0$--level set  through the saddle $(\GG, {\rm g})=(0,0)$, which in the previous section, was denoted $\cS_0$, exists only for $\d\in (0,2)$ and has equation, by~\equ{G0pl} and~\equ{delta},
\beqano\cS_0:\qquad \frac{\GG^2}{\L^2}+\frac{\rr}{a} \sqrt{1-\frac{\G^2}{\L^2}}\cos{\rm g}=\frac{\rr}{a}\ ,\eeqano
as the left hand side takes the value $\rr/a$ at the saddle. A first observation is that, if we solve for $\rr$, we find
\beqano
\rr=a\frac{\frac{\GG^2}{\L^2}}{1- \sqrt{1-\frac{\G^2}{\L^2}}\cos{\rm g}}=\frac{a(1-\ee^2)}{1-\ee\cos{\rm g}}
\eeqano
where $\ee$ is the eccentricity. This equation tells us that, while $\GG$ and ${\rm g}$ move on $\cS_0$, $\mathbf x'$ belongs to the instantaneous ellipse through $\mathbf x$, in correspondence of the true anomaly ${\rm g}$. For this reason, it is to be remarked that, while motions of $\cS_0$ are perfectly meaningful  for $\EE_0$, they might not exist for $\UU$. We also refer  $\cS_0$ as the ``collisional manifold''.

The second aspect we aim to point out  is that  the motions of $(\GG, {\rm g})$ under $\EE_0$ along this manifold can be explicitly computed.
 Indeed,
 with
\beq{conditions}\s^2:=\d(2-\d)\ ,\quad \b^2:=2-\d\ ,\quad \d\in (0,2)\eeq
  the Hamilton equation for $\GG$ is
\beqano \dst\dot\GG=\mm^2\cM\rr \sqrt{1-\frac{\GG^2}{\L^2}}\sin{\rm g}\eeqano
Eliminating ${\rm g}$ through~\equ{G0pl}, i.e.,
\beqano \sin^2{\rm g}=1-\cos^2{\rm g}=\frac{\frac{\GG^2}{\L^2}\big(\s^2-\frac{\GG^2}{\L^2}\big)}{\d^2\big(1-\frac{\GG^2}{\L^2}\big)}\eeqano
we immediately obtain the closed equation
\beqano\dot\GG=-\GG\sqrt{\L^2\s^2-\GG^2}\ .\eeqano
It can be solved by separation. The overall solution is
\beqano
\arr{
\GG(t)=\frac{\s\L}{\cosh\s\L(t-t_0)}\\
\dst{\rm g}(t)=\pm\cos^{-1}\frac{1-\frac{\b^2}{\cosh^2\s\L(t-t_0)}}{\sqrt{1-\frac{\s^2}{\cosh^2\s\L(t-t_0)}}}\ .
}\eeqano
We remark that in the case $\d=1$, hence, $\s=\b=1$, it reduces to the solution of the classical pendulum.

\subsection{Asymptotic  action--angle coordinates}\label{Asymptotic  action--angle coordinates}

The explicit construction of the action--angle coordinates for any value of $\rr$  and $\Theta$ requires the use of elliptic integrals. This is true even in the case $\Theta=0$, in which the phase portrait is, as discussed, explicit.  One possibility to deal with this situation in practical problems is to start with ``approximate coordinates''; i.e.,  coordinates such that $\EE_0$, even though being integrable, looks also as close--to--be--integrable. This would allow us to apply the machinery of perturbative methods. The first candidate to look at is the distance $\rr$ of the two fixed centers. If $\rr$ is small, $\EE_0$ is very close to $\GG^2$ and the approximate action--angle coordinates coincide with the initial $\cK$--coordinates. The case of  large $\rr$ is investigated in the present section.

When $\rr$ is large, the leading part of $\EE_0$ is
\beqano\EE_{0, 1}= \rr\mm^2{\mathcal M} \sqrt{1-\frac{\GG^2}{\L^2}} \cos{\rm g}\ .\eeqano
The integration of $\EE_{0, 1}$ relies on
solving equation
\beqa{tildecE} \sqrt{1-\frac{\GG^2}{\L^2}}\cos{\rm g}=\cE\eeqa
 for $\GG$.
The level curves~\equ{tildecE} exist only for\footnote{Observe that condition $|\cE|\le 1$ corresponds to $|\widehat\cE|\le \d$ in the notations of Section~\ref{phase portrait}.
Equation~\equ{reduced levels}  represents, in the plane $(\GG, {\rm g})$, closed curves encircling the equilibria $(0, 0)$, $(0, \p)$   (for $\cE> 0$, $\cE< 0$, respectively) corresponding to the limiting curves of the third panel in Figure~\ref{figure1} as $\d\to\infty$. } $|\cE|\le 1$ and are all closed. The associated action coordinate $\cG$ corresponds to be the area of the figure enveloped by such curves; the angle $\g$ is the time needed to reach, on a fixed level set, a given point, starting from a fixed one.  The computation of $\cG$ and $\g$ is explicit, as now we show. Solving~\equ{tildecE} for $\GG$,
 we obtain
\beqa{reduced levels}\GG=\pm\cL\sqrt{1-\frac{\cE^2}{\cos^2{\rm g}}}\quad {\rm with}\quad  |\cE|\le |\cos{\rm g}|\le 1\ ,\quad {\rm sign}(\cE)={\rm sign}(\cos{\rm g})\eeqa with $\cL=\L$.
We define the action coordinate $\cG(\cL, \cE)$ so that $\cG=0$ for $\cE=0$. Then
 \beqano
 \cG(\cL, \cE)&=&\arr{
\dst-\cL+\frac{\cL}{\p}\int_{-\arccos|\cE|}^{\arccos|\cE|}\sqrt{1-\frac{\cE^2}{\cos^2{\rm g}}}d{\rm g}\quad -1< \cE< 0\\\\
\dst\cL-\frac{\cL}{\p}\int_{-\arccos\cE}^{\arccos\cE}\sqrt{1-\frac{\cE^2}{\cos^2{\rm g}}}d{\rm g}\quad 0< \cE< 1
}
\eeqano
The period of the orbit is given  by
 \beqano\cT(\cL, \cE)=2\p \cG_{\cE}(\cL, \cE)\eeqano
With the change of variable \beqa{w}w=\frac{|\cE|}{\sqrt{1-\cE^2}}\tg{\rm g}\eeqa we obtain
\beqano
\cT(\cL, \cE)&=&4\cL|\cE|\int_{0}^{\arccos|\cE|}\frac{1}{\cos^2{\rm g}}\frac{d{\rm g}}{\sqrt{1-\frac{\cE^2}{\cos^2{\rm g}}}}=4\cL\int_{0}^{1}\frac{dw}{\sqrt{1-w^2}}=2\p\cL
\eeqano
whence (using  that $ \cG$ takes the value $0$ at $\cE=0$) the action--angle coordinates are found to be
\beqa{goodtransf} \cG=\cL\cE
\ ,\quad \g=\frac{\t}{\cL}\eeqa
where $\t$ is the time the flows employs to reach the value $(\GG, {\rm g})$ on the level set $\cE$, starting from $(\cL\sqrt{1-\cE^2}, 0)$ ($(\cL\sqrt{1-\cE^2}, \p)$, respectively).
Looking at the generating function
\beqano S(\cL, \cG, \ell, {\rm g})=\cL\ell+\int^{\rm g}\sqrt{\cL^2-\frac{\cG^2}{\cos^2{\rm g}'}}d{\rm g}'\eeqano
one obtains the  transformation of coordinates
\beqa{Gg}
\arr{\dst\L=\cL\\\\
\dst\ell=\l+\arg\left(\cos\g, \ \frac{\cL}{|\cG|}\sin\g\right)
}\quad \arr{
\dst\GG=\sqrt{\cL^2-\cG^2}\cos\g\\\\
\dst \tg{\rm g}=-\frac{\cL}{\cG}\sqrt{1-\frac{\cG^2}{\cL^2}}\sin\g\\\\  \sign\cos{\rm g}=\sign\cG}
\eeqa
The expression of $\cE$ in in action--angle coordinates is
\beqa{cE}\cE=\frac{\cG}{\cL}\ .\eeqa
Using these ``approximate'' coordinates,  one obtains the expression of $\EE_0$ as a close--to--be--integrable system for large $\rr$:
\beqano
\EE_0=\rr\mm^2\cM\frac{\cG}{\cL}+(\cL^2-\cG^2)\cos^2\g
\eeqano

\section{An application to the three--body problem (sketch) }\label{An application to the three--body problem (sketch)}
In this section we propose an application to the classical three--body problem. As said in the introduction, we choose to focus on the full problem, as opposed to the restricted one.

When the Newtonian part (i.e., the third term) in~\equ{newH2C} is much smaller compared to the Keplerian terms, one expects, by perturbation theory, that the relevant part of the dynamics of $\JJ$ is played by the ``secular terms''
\beqa{Js}\JJ_{\rm s}(\rr, \L, \Theta, \GG, {\rm g})=-\frac{\mm^3{\mathcal M}^2}{2\L^2}+ \UU(\rr, \L, \Theta, \GG, {\rm g})\eeqa
where $\UU$ is the average of the Newtonian potential, defined  in~\equ{h1}. By the previous sections, under such perturbative assumption, the motions of $\JJ$ reduce substantially to the motions of $\EE_0$.

 A  natural question is now whether a Hamiltonian  sufficiently ``close'' to $\JJ$ may generate motions that can be regarded as a ``continuation'' of the motions of  $\EE_0$.
Concretely, one might ask whether, in the specific case of the planar case,  the motions represented in Figure~\ref{figure1} may be continued to some physical system close to the planar two--centre problem. The first thought goes of course to the three--body problem.  As we are going to describe, the question demands several technical difficulties. Therefore it is just sketchily treated here.

To fix the ideas, consider the three--body problem Hamiltonian with equal masses
$$m_0=m_1=m_2$$
with the translational symmetry reduced via the heliocentric method.
Denoting as
$$\mm'=\mm=\frac{m_0}{2}\ ,\quad \cM'=\cM=2 m_0$$
the ``reduced masses'', the Hamiltonian of the system is
\beqa{3BP}
{\rm H}(  {\underline{{\mathbf y}}}, {\underline{{\mathbf x}}})&=&
\frac{\|  {{{\mathbf y}}}\|^2}{2 \mm}-\frac{\mm{\mathcal M}}{\|{{{\mathbf x}}}\|}-\frac{ \mm{\mathcal M}}{\|{{{\mathbf x}}}-{{{\mathbf x}}'}\|}-\frac{\mm'{\mathcal M}'}{\|{{{\mathbf x}}'}\|}
\nonumber\\
&+&\frac{\| {{{\mathbf y}}'}\|^2}{2\mm'}+\frac{{{{\mathbf y}}'}\cdot  {{{\mathbf y}}}}{m_0} \ .
\eeqa
We consider the Hamiltonian~\equ{3BP}, in the planar case, written in $\cK$--coordinates. Choosing the two angular momenta parallel one to the other (see Appendix~\ref{The planar case}), its expression is the following:
 \beqa{3BP1}
 \HH(\RR, \L, \GG, \rr, \ell, {\rm g})&=&-\frac{\mm^3{\mathcal M}^2}{2\L^2}-\frac{\mm{\mathcal M}}{\sqrt{{\rr }^2+2\rr  a{\rm p}+{a}^2\varrho^2}}-\frac{\mm'\cM'}{\rr}\nonumber\\
&&+
\frac{\RR^2}{2\mm'}+\frac{(\CC-\GG)^2}{2\mm'\rr^2}\nonumber\\
&&+\su{m_0}\left(\frac{\CC-\GG}{\rr}y_1(\L, \GG, \ell, {\rm g})-\RR
y_2(\L, \GG, \ell, {\rm g})
\right)
\eeqa
where $y_1(\L, \GG, \ell, {\rm g})$, $y_2(\L, \GG, \ell, {\rm g})$ are the expressions of the coordinates of the vector $\ovl{\mathbf y}$ defined in Equation~\equ{xy} below in terms of $\cK$--coordinates, given by
\beqano
y_1(\L, \GG, \ell, {\rm g})&=&\frac{\mm^2\cM}{\L(1-\ee(\L, \GG)\cos\xi(\L, \GG, \ell))}  \Big(-\sin{\rm g} \sin\xi(\L, \GG, \ell)\nonumber\\
&&+ \frac{\GG}{\L} \cos{\rm g} \cos\xi(\L, \GG, \ell)\Big)\nonumber\\
y_2(\L, \GG, \ell, {\rm g})&=&\frac{\mm^2\cM}{\L(1-\ee(\L, \GG)\cos\xi(\L, \GG, \ell))}  \Big(\cos{\rm g} \sin\xi(\L, \GG, \ell)\nonumber\\
&& + \frac{\GG}{\L} \sin{\rm g} \cos\xi(\L, \GG, \ell)\Big)
\eeqano
and  the  remaining symbols as in~\equ{Kepler Equation}--\equ{p}.

  As mentioned,   we defer a rigorous analysis of the Hamiltonian~\equ{3BP1}  to a  forthcoming paper. Here, we limit to report the results of a  numerical experiment.

The experiment has been conducted on the Hamiltonian~\equ{3BP1}.
The initial datum has been chosen as follows:
\beqa{initial datum}
&&\RR=7.071067E-005\ ,\quad   \L= 2.236067E-002\ ,\quad   \GG= 1.596860E-002\nonumber\\
&&\rr=100.000000\ ,\quad        \ell= 0.751906\ ,\quad         {\rm g}=\p\eeqa
and the total angular momentum's length $\CC=7.087036$.
The projections of the motion in the planes $({\rm g}, \GG)$, $(\ell, \L)$, $(\rr, \RR)$ are reported in Figures~\ref{figure4},~\ref{figure5} and~\ref{figure6}. We interprete  the motion of the couple $(\GG, {\rm g})$ in Figure~\ref{figure4} as a continuation of  the   motions of this couple in the last panel of Figure~\ref{figure1}, by the following heuristic considerations.
We rewrite the Hamiltonian~\equ{3BP1} as
\beqa{3BP2}
 \HH(\RR, \L, \GG, \rr, \ell, {\rm g})=\JJ(\L, \GG, \rr, \ell, {\rm g})+\KK(\RR, \rr)+f(\RR, \L, \GG, \rr, \ell, {\rm g})\eeqa
with $\JJ$ as in~\equ{EEJJ} with $\Theta=0$ and
\beqa{fullH}
\KK(\RR, \rr)&:=&\frac{\RR^2}{2\mm'}+\frac{\CC^2}{2\mm'\rr^2}-\frac{\mm'\cM'}{\rr}\nonumber\\
 f(\RR, \L, \GG, \rr, \ell, {\rm g})&:=&\frac{-2\CC\GG+\GG^2}{2\mm\rr^2}+\su{m_0}\Big(\frac{\CC-\GG}{\rr}y_1(\L, \GG, \ell, {\rm g})\nonumber\\
&&-\RR
y_2(\L, \GG, \ell, {\rm g})
\Big)
\eeqa
The large gap between the initial size of $a=10^{-3}$ in~\equ{p}  and $\rr=10^2$ makes the full Hamiltonian~\equ{3BP2} to be well represented by its $\ell$--average
\beqano
\ovl\HH(\RR, \L, \GG, \rr, \ell, {\rm g})&=&-\frac{\mm^2\cM^2}{2\L^2}+\UU(\L, \GG, \rr, {\rm g})+\KK(\RR, \rr)+\frac{-2\CC\GG+\GG^2}{2\mm\rr^2}\eeqano
having used the last term in $f$ in~\equ{fullH} has zero average.
When the last term is neglected, $\ovl\HH$ reduces to
\beqano\hh(\L, \GG, \rr, {\rm g})=-\frac{\mm^2\cM^2}{2\L^2}+\UU(\L, \GG, \rr, {\rm g})+\KK(\RR, \rr)\ .\eeqano
This Hamiltonian is not integrable. However, $\KK(\RR, \rr)$ has an equilibrium at $(\RR, \rr)=\left(0, \frac{\CC^2}{{\mm'}^2\cM'}\right)$. If $\CC$ is sufficiently large, such equilibrium  is an ``approximate'' equilibrium to $\hh$, so the motions of $\hh$ are approximately decoupled, and are obtained combining small oscillations (generated by $\KK$) of the couple  $(\RR, \rr)$ about the equilibrium with the motions (generated by $\EE_0$) for $(\GG, {\rm g})$ depicted in the last panel of Figure~\ref{figure1} (since $\d= 10^{5}$), with Keplerian motions for the couple $(\L, \ell)$ generated by the Keplerian part in~\equ{Js}.  Observe that the initial values of $\RR$ and $\rr$ in~\equ{initial datum}  are very close to the ones at the equilibrium (with the above choices of $\CC$ and of the masses, $\frac{\CC^2}{{\mm'}^2\cM'}=100.452159$), so we expect that this picture of the motion is preserved in $\HH$ for a long time, at least on a positive measure set of initial data. However, a rigorous treatment of these arguments is a bit delicate for the following reasons. A first difficulty is that the unperturbed motions of $\GG$ are very close to $\GG=0$. This occurrence makes the Hamiltonian singular, even though the unperturbed part is not so. At this singularity we ascribe the divergences of $\L$ in Figure~\ref{figure5}. Another difficulty is that, as   $\rr$ is large compared to $a$, besides the perturbative term, $\UU$ is also small. Therefore,  in order  that the motions of $(\GG, {\rm g})$ are in fact ``led'' by $\UU$, a careful balance between $\UU$ and $\ovl f:=\frac{-2\CC\GG+\GG^2}{2\mm\rr^2}$ is required. We should not miss to mention another delicate question. When $0<\d<2$ (the first and the second panel in Figure~\ref{figure1}),  the previous heuristic arguments do not seem to  apply. Moreover, as outlined in Section~\ref{The collisional manifold and its motions}, the separatrix $\cS_0$ appearing in such figures is a singular manifold for $\UU$. These considerations convinced the author that a rigorous treatment of  a possible application of the results of Sections~\ref{coordinates}--\ref{Dynamical consequences} to the three--body problem  would require such an accurate analysis to exceed  the purposes of the present note.
\section{Conclusion}
\begin{figure}
\vspace*{10pt}\center{ \includegraphics[height=5.0cm, width=5.0cm]{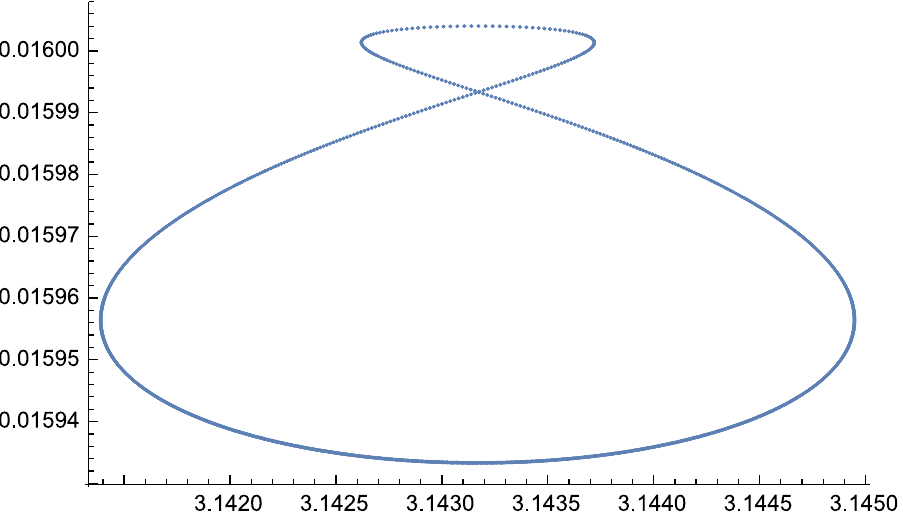}}  %
\caption{Projection of the motion in the plane (g, G). }\label{figure4}
\vspace*{10pt} \center{\includegraphics[height=5.0cm, width=5.5cm]{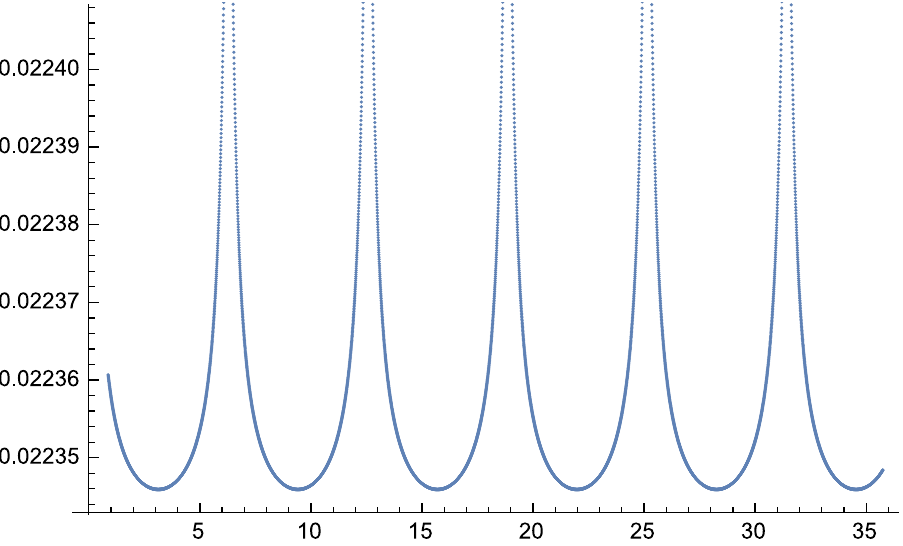}}  \caption{
Projection of the motion in the plane $(\ell, \L)$.}\label{figure5}
\vspace*{10pt}\center{ \includegraphics[height=5.0cm, width=9.0cm]{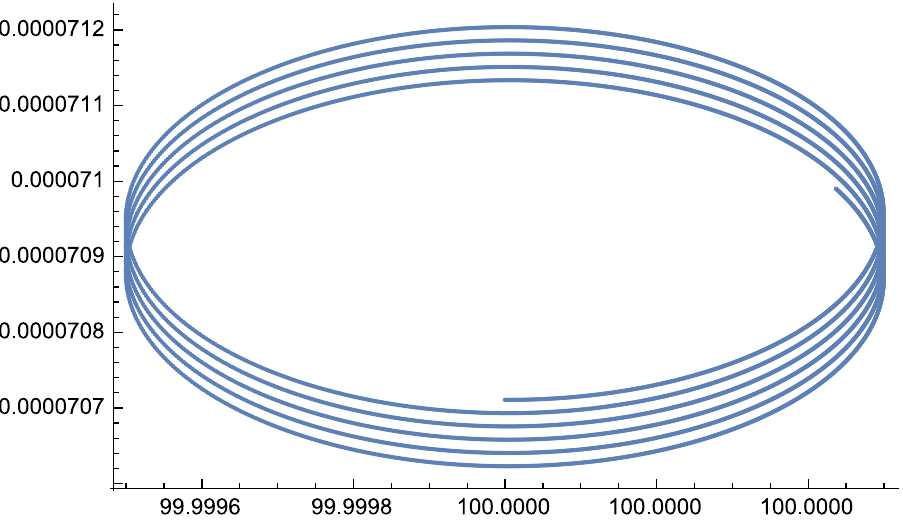}}
 \caption{Projection of the motion in the plane $(\rr, \RR)$. }\label{figure6}
 \end{figure}

We propose an alternative approach to the analysis of the two--centre problem. We dismiss the separability property~\cite{bekovO78} that one usually gains using the ellipsoidal coordinates. Rather, we regard the two--centre Hamiltonian $\JJ$ in~\equ{newH2C}  as a small perturbation of the Kepler problem. This is possible when the primaries are very far or their masses   are sensitively  different. We analyse the Hamiltonian of the problem using a special set of canonical coordinates, denoted as $\cK$, and defined in Section~\ref{coordinates}, in terms of which the Hamiltonian has two degrees of freedom. This is because $\cK$ includes, among its coordinates, all the first integrals of $\JJ$ but one. The lack of separation is compensated by the fact that, as proved in~\cite{pinzari19}, the associated secular Hamiltonian in {\it renormalizably integrable}, meaning  that it can be expressed as a function of a ``normalising'', much simpler, function. The phase portrait of such normalising function can be studied exactly, at least in the case of the planar problem, in correspondence of any value of the ratio $\d=\rr/a$, where $\rr=\|\bx'\|$ and $a$ is the semi--major axis of the instantaneous ellipse generated by the Keplerian part of $\JJ$ (Section~\ref{phase portrait}; Figure~\ref{figure1}). Such phase portrait shows, for any value of $\d$, the existence of small oscillations of the perihelion of the instantaneous ellipse accompanied by a periodic change of the shape of the ellipse that have large eccentricity and, at every period, squeezes to a segment while the body changes its direction on it. We call such motions {\it perihelion librations}. Moreover, when $0<\d<2$, the phase portrait includes  a saddle equilibrium point and a separatrix $\cS_0$ through it.
 It is also possible to compute, exactly, the motion on $\cS_0$  (Section~\ref{The collisional manifold and its motions}) as well as at least a first order approximation of the  action--angle coordinates in the case that $\d$ is very large (Section~\ref{Asymptotic  action--angle coordinates}). In Section~\ref{An application to the three--body problem (sketch)} we conjecture it is possible to use the perturbative approach proposed in the paper in order  to prove, in the classical three--body problem, the existence motions including perihelion librations, at least in the case that $\d$ is very large. We propose one numerical experiment and a heuristic argument that seem to evidence the conjecture (Figures~\ref{figure4},~\ref{figure5} and~\ref{figure6}).

\appendix

 \section{Explicit formulae of the ${\cK}$--map}\label{explicit formulae}
 Here we provide the analytical expression of the map
 \beqa{K map}
 \phi:\quad {{\cK}}=\Big(\ZZ,\CC, \Theta, \GG, \L, \RR  , \zeta, g, \vartheta, {\rm g}, \ell, \rr \Big)
 \to (\underline{\mathbf y}, \underline{\mathbf x})=({\mathbf y}', {\mathbf y}, {\mathbf x}', {\mathbf x})\ .
 \eeqa
 which is needed to write $\JJ$ and $\EE$ in  terms of ${\cK}$.

Let
\beqa{k incli} i=\cos^{-1}\left(\frac{\ZZ}{\CC}\right)\ ,\quad  i_1=\cos^{-1}\left(\frac{\Theta}{\CC}\right)\ ,\quad i_2=\cos^{-1}\left(\frac{\Theta}{\GG}\right)\eeqa
and
{\small\beqa{R1R3}
&&\cR_1(\a):=\left(
\begin{array}{ccc}
1&0&0\\
0&\cos \a&-\sin\a\\
0&\sin\a&\cos\a
\end{array}
\right)\quad  \cR_3(\a):=\left(
\begin{array}{ccc}
\cos\a&-\sin\a&0\\
\sin\a&\cos\a&0\\
0&0&1
\end{array}
\right)\ .\eeqa}Using the definitions in~\equ{belline} and the observation that, if $\FF\to^{(\YY, \XX,  x)}\FF'$, the transformation of coordinates
which relates the coordinates ${\mathbf x}'$ relatively to $\FF'$ to  the coordinates ${\mathbf x}$ relatively to $\FF$  is
\beqano {\mathbf x}=\cR_3(x)\cR_1(\iota){\mathbf x}'\eeqano
where $\cR_1$, $\cR_3$ are as in~\equ{R1R3}, while $\iota:=\cos^{-1}\frac{\YY}{\XX}$, for the map~\equ{K map} we find the following analytical expression
 \beqa{yx}
\phi:\qquad \left\{
\begin{array}{ll}
{\mathbf x}= \cR_3(\zeta)\cR_1(i)\cR_3(g)\cR_1(i_1)\cR_3(\vartheta)\cR_1(i_2)\ovl{\mathbf x}(\L,\GG,\ell, {\rm g})\\ \\
{\mathbf y}= \cR_3(\zeta)\cR_1(i)\cR_3(g)\cR_1(i_1)\cR_3(\vartheta)\cR_1(i_2)\ovl{\mathbf y}(\L,\GG,\ell, {\rm g})
\\\\
{\mathbf x}'=\rr  \cR_3(\zeta)\cR_1(i)\cR_3(g)\cR_1(i_1){\mathbf k}\\ \\
 {\mathbf y}'=\frac{{\rm R}'}{{\rm r}'}{\mathbf x}'+\frac{1}{{\rm r}'^2}{\mathbf M}'\times{\mathbf x}'
\end{array}
\right.
\eeqa
with
\beqa{xy}
&&
\ovl{\mathbf x}=
\frac{\L^2}{\mm^2{\mathcal M}}\cR_3({\rm g}-\p/2) \left(
\begin{array}{ccc}
\cos\xi(\L,\GG,\ell)-\ee\\
\frac{\GG}{\L}\sin\xi(\L,\GG,\ell)\\
0
\end{array}
\right)\nonumber\\
&& \ovl{\mathbf y}=\frac{\mm^2{\mathcal M}}{\L(1-\ee(\L, \GG)\cos\xi(\L, \GG, \ell))}\cR_3({\rm g}-\p/2) \left(
\begin{array}{ccc}
-\sin\xi(\L,\GG,\ell)\\
\frac{\GG}{\L}\cos\xi(\L,\GG,\ell)\\
0
\end{array}
\right)\nonumber\\
&&{\mathbf C}=\CC \cR_3(\zeta)\cR_1(i){\mathbf k}\ ,\quad
{\mathbf M}=\GG \cR_3(\zeta)\cR_1(i)\cR_3(g)\cR_1(i_1)\cR_3(\vartheta)\cR_1(i_2){\mathbf k}\nonumber\\
&&{\mathbf M}'={\mathbf C}-{\mathbf M}\eeqa
where  $\ee(\L, \GG)$
and $\xi(\L, \GG, \ell)$ are as in~\equ{Kepler Equation}--\equ{p}. 

\subsection{The planar case}\label{The planar case} The planar motions are obtained setting $\Theta=0$ and $\vartheta=0$, $\p$. Indeed, the planar motions correspond to take ${\mathbf C}\parallel\s {\mathbf M}=\s({\mathbf C}-{\mathbf M}')$, with $\s=\pm 1$, so
\beqano\Theta=\frac{{\mathbf M}\cdot{\mathbf x}'}{\|{\mathbf x}'\|}=\mathbf 0\ .\eeqano Moreover, from the definitions in~\equ{nodes}, we have ${\mathbf i}_2\parallel (-\s{\mathbf i}_3)$, so \beqano\vartheta=\a_{{\mathbf x'}}({\mathbf i}_2, {\mathbf i}_3)=\arr{\dst \p\quad \rm if\quad \s=+1\\\\
\dst 0\quad \rm if\quad \s=-1\ .
}\eeqano
In such  cases, the map~\equ{K map} reduces to
\beqa{phipl}
 \phi_{\rm pl}:\quad {{\cK}}=\Big(\ZZ,\CC, \GG, \L, \RR  , \zeta, g,  {\rm g}, \ell, \rr \Big)
 \to (\underline{\mathbf y}, \underline{\mathbf x})=({\mathbf y}', {\mathbf y}, {\mathbf x}', {\mathbf x})\ .
 \eeqa
with
\beqa{yy'}
\left\{
\begin{array}{lll}
\dst {\mathbf x}=\cR_3(\zeta)\cR_1(i)\cR_3(g)(\ovl{\mathbf x}_1{\mathbf i}+\ovl{\mathbf x}_2{\mathbf j})\\\\
\dst {\mathbf y}=\cR_3(\zeta)\cR_1(i)\cR_3(g)(\ovl{\mathbf y}_1{\mathbf i}+\ovl{\mathbf y}_2{\mathbf j})\\\\
\dst {\mathbf x}'=-\rr\cR_3(\zeta)\cR_1(i)\cR_3(g){\mathbf j}\\\\
\dst {\mathbf y}'=-\RR\cR_3(\zeta)\cR_1(i)\cR_3(g){\mathbf j}+\frac{\CC-\s\GG}{\rr}\cR_3(\zeta)\cR_1(i)\cR_3(g){\mathbf i}
\end{array}
\right.
\eeqa
We call {\it planar and prograde} the case $\s=+1$; while {\it planar and retrograde} the case $\s=-1$.

 \subsection{Derivation of the formulae~\equ{EEJJ}}\label{JEK}

 Using the general formulae in~\equ{yx}, we find
\beqa{inner0}{\mathbf x}'\cdot {\mathbf x}
&=&{\mathbf k}\cdot \cR_1(i_2)\ovl{\mathbf x}(\L,\GG,\ell, {\rm g})=-\rr a\sqrt{1-\frac{\Theta^2}{\GG^2}}{\rm p}\eeqa
with ${\rm p}$ as in~\equ{p},
and where we have used $\cR_3^{\rm t}(\vartheta){\mathbf k}={\mathbf k}$, the relation
\beqano\sin i_2=\sqrt{1-\frac{\Theta^2}{\GG^2}}\eeqano
(which is implied by the definition of $i_2$ in~\equ{k incli}) and the expression for $\ovl{\mathbf x}$ in~\equ{xy}.
Equations~\equ{inner0},~\equ{xy} and the definition of $\rr =\|{\mathbf x}'\|$
 then imply that
the Euclidean distance between ${\mathbf x}'$ and ${\mathbf x}$ has the expression
\beqa{distance}\|{\mathbf x}'-{\mathbf x}\|^2={{\rr }^2+2\rr  a\sqrt{1-\frac{\Theta^2}{\GG^2}}
{\rm p}+a^2\varrho^2}\ .\eeqa
Recall that the eccentricity vector ${\mathbf L}$ in~\equ{CL} is related to $\ee$ and ${\mathbf P}$ via  ${\mathbf L}=\mm^2{\mathcal M}\ee{\mathbf P}$. The expression of ${\mathbf P}$ is obtained from the one for ${\mathbf x}$
in~\equ{yx} taking $\n=\ell=0$ and normalizing. Namely,
\beqano{\mathbf P}=\frac{\ovl{\mathbf x}}{a\varrho}=\cR_3(\zeta)\cR_1(i)\cR_3(g)\cR_1(i_1)\cR_3(\vartheta)\cR_1(i_2)\ovl{\mathbf P}
\eeqano
with \beqano\ovl{\mathbf P}=\left(
\begin{array}{ccc}
\sin{\rm g}\\
-\cos{\rm g}\\
0
\end{array}
\right)\eeqano
Then, analogously to~\equ{inner0},
 we find, for the inner product ${\mathbf x}'\cdot {\mathbf P}$ the expression
\beqa{inner}{\mathbf x}'\cdot {\mathbf P}=-\rr \sqrt{1-\frac{\Theta^2}{\GG^2}}\cos{\rm g}\eeqa
Using the formulae in~\equ{distance},~\equ{inner}, the definition of $\GG=\|\mathbf x\times \mathbf y\|$,
 we find that
the functions ${\JJ}$, $\EE$, written terms of the coordinates ${{\cK}}$, are as in~\equ{EEJJ}.

\section*{Acknowledgments} I heartily thank Krzysztof Cieplinski and Maciej Capinski for their interest and the anonymous Reviewers for  useful suggestions. Thanks also to Alessandra Celletti,  Amadeu Delshams and Susanna Terracini for their interest;  Marcel Guardia and
Tere Seara for fruitful discussions. I am indebted to Jer\^ome Daquin, who carefully read the code I used to produce  Figures~\ref{figure4},~\ref{figure5} and~\ref{figure6}, and pointed  out a misprint.
Figures~\ref{figure1}--\ref{figure6} have been drawn with {\sc mathematica}\textsuperscript{\textregistered}.

\end{document}